\documentclass[11pt,reqno]{article}

\usepackage{amsfonts}
\usepackage{amsmath}
\usepackage{amsthm}
\usepackage{amssymb}
\usepackage{mathrsfs}
\usepackage{amstext}

\usepackage[dvips]{graphicx}

\font\msbm=msbm10

\numberwithin{equation}{section}

\theoremstyle{plain}
\newtheorem{Theorem}{Theorem}[section]
\newtheorem{lemma}[Theorem]{Lemma}
\newtheorem{corollary}[Theorem]{Corollary}

\newtheorem{remark}[Theorem]{Remark}
\def\mathbb#1{\hbox{\msbm{#1}}}
\newcommand{\N}{{\mathbb{N}}}

\newcommand{\Z}{{\mathbb{Z}}}
\newcommand{\C}{{\mathbb{C}}}

\renewcommand{\P}{{\mathbb{P}}}
\newcommand{\E}{{\mathbb{E}}}
\newcommand{\B}{{\cal{B}}}
\newcommand{\A}{{\cal{A}}}

\newcommand{\sgn}{\operatorname{sgn}}
\newcommand{\Tr}{\operatorname{Tr}}

\newcommand{\beq}{\begin{eqnarray}}
\newcommand{\eeq}{\end{eqnarray}}
\newcommand{\beqn}{\begin{eqnarray*}}
\newcommand{\eeqn}{\end{eqnarray*}}

\newcommand{\F}{{\cal F}}
\newcommand{\ol}{\overline}
\newcommand{\supp}{\operatorname{supp}}

\renewcommand{\Re}{\operatorname{Re}}

\renewcommand{\qed}{\rule{2.5mm}{2.5mm}}
\newenvironment{Proof}{\noindent
{\bf\underline{Proof:} }}
{\hspace*{\fill}\qed\vskip1em}
\begin{document}
\title{Random Sampling of Sparse Trigonometric Polynomials}

\author{Holger Rauhut\\
NuHAG, Faculty of Mathematics, University of Vienna\\
Nordbergstrasse 15, A-1090 Wien, Austria\\
rauhut@ma.tum.de}

\maketitle

\begin{abstract}
We study the problem of reconstructing a multivariate 
trigonometric polynomial having only few non-zero coefficients 
from few random samples.
Inspired by recent work of
Candes, Romberg and Tao we propose 
to recover the polynomial by Basis Pursuit, i.e., by $\ell^1$-minimization. 
Numerical experiments show that in many cases the trigonometric polynomial
can be recovered exactly provided
the number $N$ of samples is high enough compared
to the ``sparsity'' -- the number of non-vanishing coefficients. However,
$N$ can be chosen small compared to the assumed maximal degree of the 
trigonometric polynomial.
Hence, the proposed scheme may overcome the Nyquist rate.
We present two theorems that explain this observation.
Unexpectly, they establish a connection to an interesting combinatorial problem
concerning set partitions, which seemingly has not yet been considered before. 
\end{abstract}
\vspace{0.5cm}



\section{Introduction}

Recently, Candes, Romberg and Tao observed the surprising fact that
it is possible to recover certain discrete functions exactly from 
vastly incomplete information on their discrete Fourier 
transform \cite{CRT1,CT,CRT2,CR}. The crucial property of these functions is
their sparsity with respect to the canonical (Dirac) basis, i.e., their 
(unknown) support is very small. 
The recovery procedure consists in minimizing the $\ell^1$-norm
of the signal subject to the constraint that the Fourier coefficients are matched. 
This task is also known as Basis Pursuit \cite{CDS}.
Since minimizing the total variation norm can be reformulated as
minimizing the $\ell^1$-norm there are relevant applications
in image processing, in particular, computer tomography \cite{CRT1,CR}.

This paper is concerned with the related problem of reconstructing
a sparse trigonometric polynomial from few randomly chosen samples drawn
from the continuous uniform distribution on $[0,2\pi]^d$.
By ``sparse'' we mean that only very few coefficients of the 
polynomial are non-zero. However, a priori we do not know the 
support of the coefficients. From a practical viewpoint considering such
polynomials can be motivated as follows. First, trigonometric polynomials
with a certain maximal degree model band-limited signals.
Secondly, in many cases it seems reasonable that only few
coefficients (with unknown location) are large. Such a signal can at least 
be approximated by a sparse one. 

We propose to reconstruct the sparse polynomial
from its random samples by Basis Pursuit similarly 
as in \cite{CRT1,CT,CRT2,CR}. 
From numerical experiments
it is evident that this scheme can indeed reconstruct the polynomial
exactly provided the number of samples is large enough with respect to
the sparsity.
When comparing the number of samples to the assumed maximal degree
of the polynomial it turns out that this method may overcome 
the Nyquist rate by far. Thus, the described recovery method is very likely to have
high potential for practical applications in signal processing.

We will present two theorems that explain the observed phenomenon.
Similar to \cite{CRT1} the first one
estimates the probability of exact reconstruction given an arbitrary
sparse trigonometric polynomial. Hence, this can be viewed as a worst case
estimate.
Our second theorem is more directed towards an avarage case analysis.
It gives a probability estimate for 
generic polynomials in the sense that the support of the 
coefficients is modelled as random set. A result of this type seems to be new.
As one may expect it gives better probability estimates than the first one.
However unexpectly, it relates the problem to a seemingly new and
difficult combinatorial
problem about set partitions. Unfortunately, we were not 
able to solve this problem in general, and as a consequence 
we cannot yet exploit fully our probability estimate.
We have to leave the combinatorial aspect as an interesting
open problem.
 
We would like to mention some recent related work.
In \cite{CT,CRT2} Candes et al. study stability aspects of the problem and
investigate also recovery from few inner products with
random vectors following Gaussian distributions 
and binary distributions.
In \cite{CR} some practical examples are presented. 
The recovery from Gaussian measurements via Basis Pursuit
is also investigated by Rudelson and Vershynin in \cite{RV} in the context of
error correcting codes, while
Tropp \cite{Tropp} studies the reconstruction by 
Orthogonal Matching Pursuit. 
In \cite{Donoho1,DTsaig} Donoho and Tsaig introduce the 
terminology ``compressed sensing'' for this range of problems 
and in \cite{Donoho2,DT} probabilistic results concerning 
Basis Pursuit are discussed.
A randomized sublinear algorithm for reconstructing sparse Fourier data is 
introduced and analyzed in \cite{ZGSD}.
If the reader is interested in reconstructing not necessarily sparse
trigonometric polynomial from random samples we refer to recent work of
Bass and Gr\"ochenig \cite{BG}, where probabilistic estimates
of related condition numbers are developed.

The paper is structured as follows. In Section \ref{sec_Main_Results}
we describe the problem and present our main results. To this end
we also need to introduce some background on set partitions.
Section \ref{sec_Proof} will be devoted to the proofs. 
Section \ref{sec_Qn} gives some more information on
the combinatorial problem related to our second theorem.
In Section \ref{sec_bound} we present some plots of the probability bounds
resulting from our theorems and finally Section \ref{sec_numerical}
describes some numerical experiments.  

\bigskip 

\noindent

{\bf Acknowledgements:} The author was supported by the 
European Union's Human Potential Programme, under contract
HPRN--CT--2002--00285 (HASSIP). He would like to thank Stefan Kunis for 
stimulating discussions on numerical aspects of the topic. 
Also he acknowledges interesting conversations with Justin Romberg and
his mail correspondence with Emmanuel Candes on the subject.

\section{Description of the Main Results}
\label{sec_Main_Results}

\subsection{The Setting}

Let $\Pi_q = \Pi_q^d$ denote the space of all trigonometric polynomials
of maximal order $q \in \N_0$ in dimension $d$. An element $f$ of $\Pi_q$
is of the form 
\[
f(x) \,=\, \sum_{k \in [-q,q]^d \cap \Z^d} c_k e^{i k\cdot x}, \qquad 
x \in [0,2\pi]^d, 
\]
with some Fourier coefficients $c_k \in \C$.
The dimension of $\Pi_q^d$ will be denoted by $D:= (2q+1)^d$. In the
sequel we will use the short notation $[-q,q]^d$ instead of 
$[-q,q]^d \cap \Z^d$.

Through the rest of this paper we will be dealing with
``sparse'' trigonometric polynomials, i.e., we assume that the sequence of 
coefficients $c_k$ is supported only on a set $T$, which is much smaller
than the dimension $D$ of $\Pi_q$. However, a priori nothing is known about
$T$ apart from a maximum size. Thus, it is useful to introduce
the set $\Pi_q(M)=\Pi_q^d(M) \subset \Pi_q$ of all trigonometric polynomials
whose Fourier coefficients are supported on a set 
$T \subset [-q,q]^d \cap \Z^d$ satisfying 
$|T| \leq M$, i.e., $f \in \Pi_{q}(M)$ is of the form 
$f(x) = \sum_{k \in T} c_k e^{ik\cdot x}$. Note that $\Pi_q(M)$ is not
a linear space.

Our aim is to sample a trigonometric polynomial $f$ of $\Pi_q(M)$ at $N$
randomly chosen points and try to reconstruct $f$ from these
samples. We model the sampling points
$x_1,\hdots,x_N$ as independent random variables having the uniform
distribution on $[0,2\pi]^d$. We collect them into the sampling set 
\[
X\,:=\,\{ x_1,\hdots,x_N \}. 
\]
Obviously, the cardinality $|X|$ equals the number of samples
$N$ with probability $1$.

Motivated by results of Candes, Romberg and Tao
in \cite{CRT1} we propose the following non-linear
method of reconstructing
$f \in \Pi_d(M)$ from its sampled values $f(x_1),\hdots,f(x_N)$. 
We minimize the $\ell_1$-norm of the Fourier coefficients $c_k$,
\[
\|(c_k)\|_{1} \,:=\, \sum_{k \in [-q,q]^d} |c_k|,
\]
under the
constraint that the corresponding trigonometric polynomial matches
$f$ on the sampling points. That is we solve the problem
\begin{equation}\label{min_problem}
\min \|(c_k)\|_{1}
\quad\mbox{ subject to }\quad
g(x_j) \,:=\, \sum_{k \in [-q,q]^d} c_k e^{ik\cdot x_j} \,=\, f(x_j), 
\quad j=1,\hdots,N.
\end{equation}
This task -- also referred to as Basis Pursuit \cite{CDS} -- 
can be performed with efficient convex optimization techniques \cite{BV},
or even linear programming in case of real-valued coefficients $c_k$.

Once all the coefficients $c_k$, $k\in [-q,q]^d$, 
are known also $f$ is known completely and 
can be evaluated efficiently at any point, e.g., with
the Fourier transform for non-equispaced data developed by Daniel Potts
et al. \cite{DPT}.

\begin{figure}\label{fig_reconstruct}
\parbox[t]{8cm}{
\includegraphics[width=8cm]{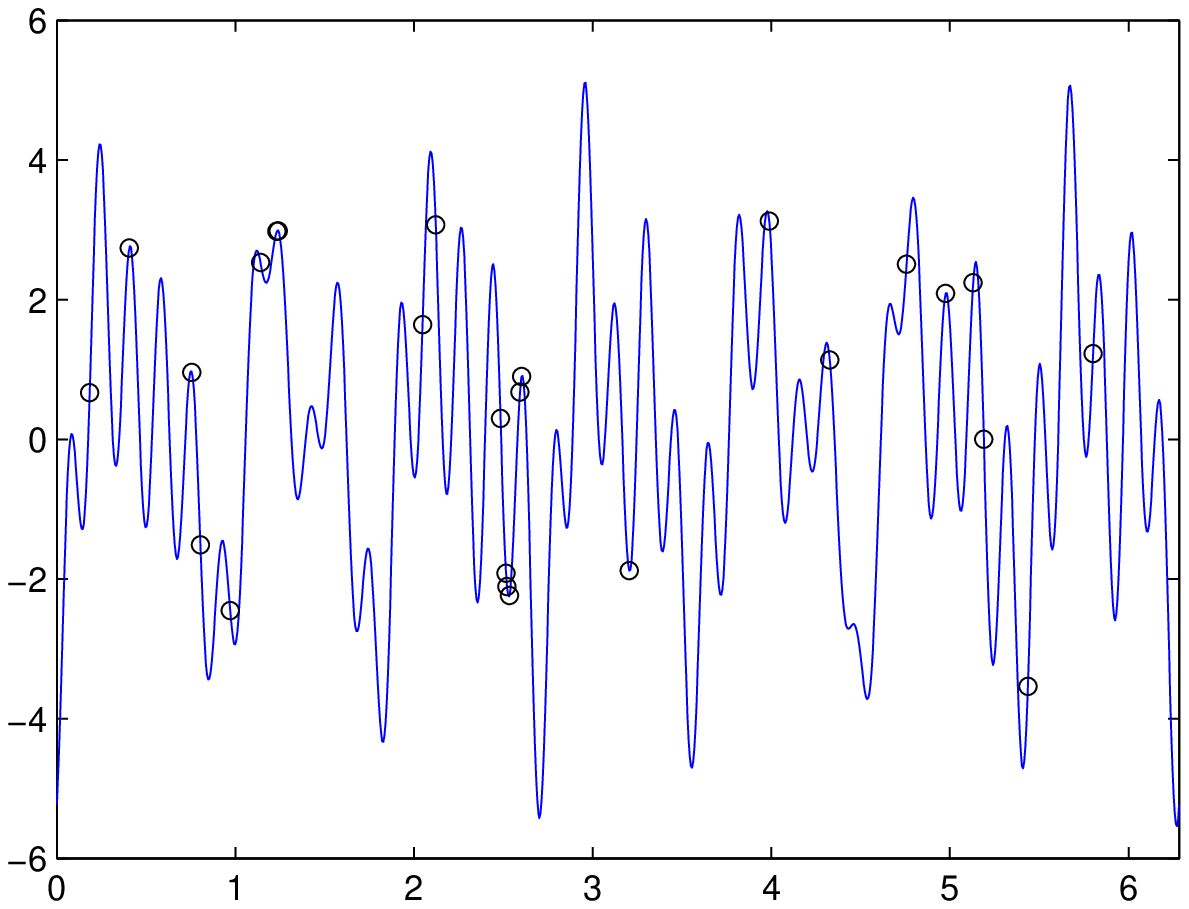}}
\parbox[t]{8cm}{
\includegraphics[width=8cm]{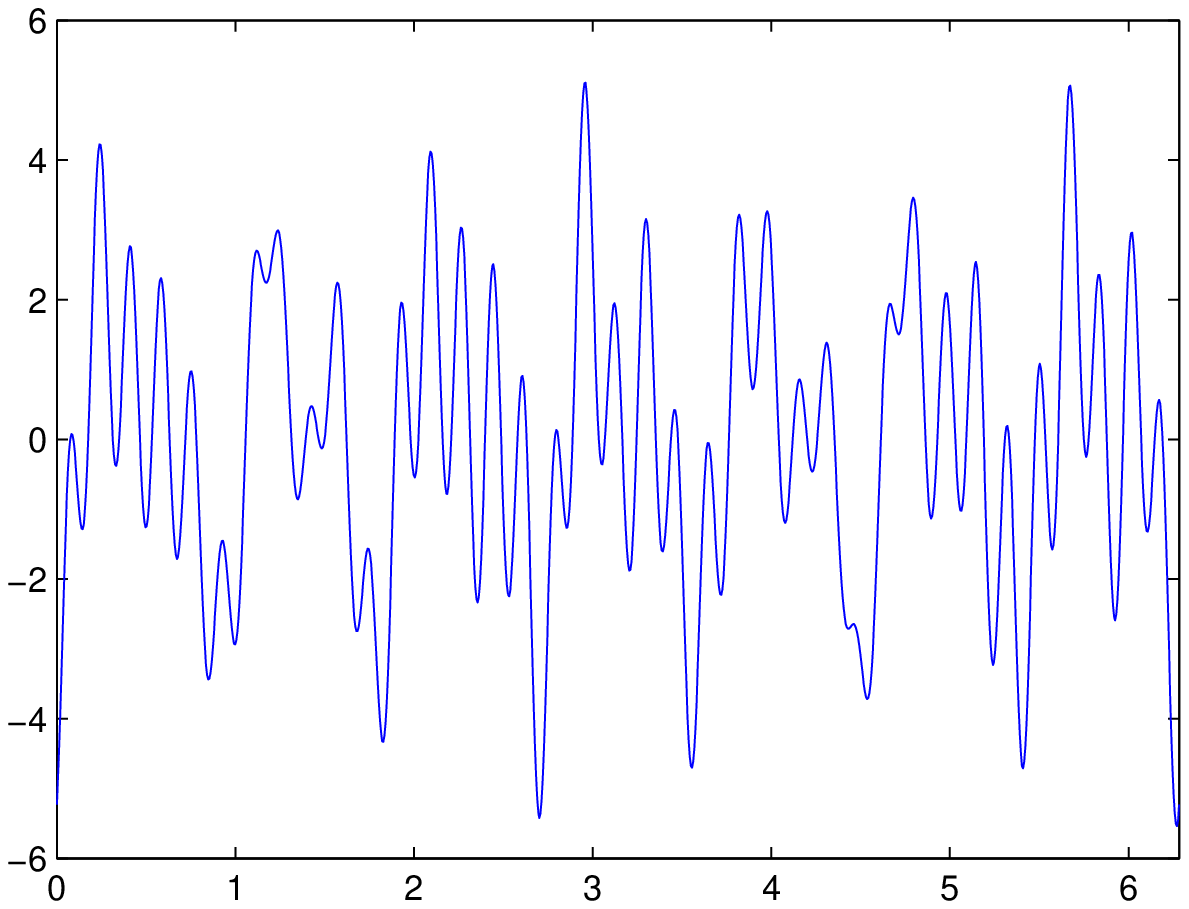}}
\caption{Original sparse trigonometric polynomial with samples (left)
and reconstruction (right)}
\end{figure}

Surprisingly, there is numerical evidence that the above reconstruction 
scheme recovers $f$ {\it exactly} provided the number 
of samples is large enough
compared to the sparsity. Indeed, Figure \ref{fig_reconstruct}
shows a sparse trigonometric polynomial with $8$ non-zero coefficients
and $N=25$ sampling points while the maximal degree is $q=40$, i.e.,
$D=81$. The right
hand side shows the reconstruction from the samples by solving the minimization
problem (\ref{min_problem}). The reconstruction is exact!
We refer to Section \ref{sec_numerical} for more information
on the numerical experiments.

Our main results are two theorems that give a theoretical explanation
of this phenomenon. The first one treats any sparse polynomial
in $\Pi_q^d(M)$ and the second one considers ``generic'' polynomials
in the sense that the set $T$ of non-vanishing coefficients
is modelled as random set.
Unexpectly, both results involve combinatorial quantities connected to set
partitions. We will spend the next section introducing the necessary
notation.
 
\subsection{Set Partitions}

We denote $[n]:= \{1,2,\hdots,n\}$. A partition of $[n]$ is 
a set of subsets of $[n]$ -- called blocks -- 
such that each $j \in [n]$ is contained
in precisely one of the subsets. By $P(n,k)$ we denote the set of all
partitions of $[n]$ into exactly $k$ blocks such that each block
contains at least $2$ elements. For example $P(4,2)$ consists
of 
\[
\{\{1,2\},\{3,4\}\},\quad \{\{1,3\},\{2,4\}\}\quad 
\mbox{and } \quad \{\{1,4\},\{2,3\}\}.
\]
Clearly, $P(n,k)$ is empty if $k > n/2$. The numbers $S_2(n,k)=|P(n,k)|$ 
are called associated Stirling numbers of the second kind. 
They have the following exponential 
generating function, see 
\cite[formula (27), p.77]{Riordan},
\begin{equation}\label{expgen_func}
\sum_{n=1}^\infty \sum_{k=1}^{\lfloor n/2 \rfloor}  S_2(n,k) y^k \frac{x^n}{n!}
\,=\, \exp(y(e^x-x-1)).
\end{equation}
Based on this one may deduce that
the numbers $S_2(n,k)$ satisfy the recursion formula
\begin{equation}\label{rec_formula}
S_2(n,k) \,=\, k S_2(n-1,k) + (n-1)S_2(n-2,k-1).
\end{equation}
Also a combinatorial argument for this recursion exists, see 
Section \ref{sec_Qn} where also further
information on the numbers $S_2(n,k)$ will be given.

We also need partitions of a different type. An adjacency is defined to be
an occurence of two consecutive integers of $[n]$ in the same block of
a partition. Hereby, consecutive is understood in the circular sense, i.e.,
also $n$ and $1$ are considered consecutive. We define $U(n,k)$
as the set of all partitions into $k$ subsets 
having {\it no} adjacencies. For instance, $U(5,3)$
consists of the partitions
\begin{align}
\{\{1,4\},\{2,5\},\{3\}\}&,\quad \{\{1,4\},\{2\},\{3,5\}\},\quad \{\{1\},\{2,4\},\{3,5\}\},\notag\\
\{\{1,3\},\{2,5\},\{4\}\}&\quad \mbox{ and }\quad \{\{1,3\},\{2,4\},\{5\}\}.
\end{align}
Clearly, $U(n,1)$ is empty.
We remark that it was only very recently 
that D. Knuth \cite{Knuth} raised the problem
of determining the number of partitions in $U(n,k)$. 

We will also need a slight variation of this type of partitions.
Let $[K] \times [m] = \{1,\hdots, K\} \times \{1,\hdots,m\}$ for some
numbers $K,m \in \N$. We denote
by $U^*(K,m,s)$ the set of all partitions of $[K] \times [m]$
such that $(p,u)$ and $(p,u+1)$ are not contained in the same block
for all $p=1,\hdots,K$ and $u=1,\hdots,m-1$. (So this sort of 
consecutiveness is {\it not} understood in the circular sense.)
We remark that $U(K,1,k)$ is the set of all partitions
of a $K$-element set into $k$ subsets (without any restriction on the type
of partition). In particular, the numbers $|U(K,1,k)|$
equal the (ordinary) Stirling numbers $S(K,k)$ of the second kind. The numbers
$b_n := \sum_{k=1}^n S(n,k)$ are called Bell numbers \cite{Riordan,Stanley}. 

Now let $\A = \{A_1,\hdots,A_t\}$
be a partition in $P(n,t)$ and $\B = \{B_1,\hdots,B_s\} \in U(n,s)$.
By $A_i + 1$ we understand the set whose elements are the ones of $A_i$ 
incremented by $1$ in the circular sense, i.e., $n+1 \equiv 1$.
We associate a $t\times s$ matrix $M=M(\A,\B)$ to the pair $\A,\B$ by
setting
\begin{equation}\label{defAB}
M_{i,j} \,:=\, |A_i \cap B_j| - |(A_i+1) \cap B_j|,
\qquad 1\leq i\leq t,\, 1\leq j\leq s.
\end{equation}
Then we define $Q(n,t,s,R)$ to be the number of 
pairs of partitions $(\A,\B)$ with
$\A \in P(n,t)$ and $\B \in U(n,s)$ such that 
the rank of $M(\A,\B)$ equals $R$, i.e.,
\begin{equation}\label{def_Q}
Q(n,t,s,R) \,:=\, \# \left\{(\A,\B):\, \A \in P(n,t), \B \in U(n,s), 
\operatorname{rank} M(\A,\B) = R \right\}.
\end{equation}
Observe that 
\[
\sum_{i = 1}^t M_{i,j} \,=\, \sum_{i=1}^t 
\left(|A_i \cap B_j| - |(A_i + 1) \cap B_j|\right)
\,=\, |\{1,\hdots,n\} \cap B_j| - |\{1,\hdots,n\} \cap B_j| \,=\, 0
\] 
(since the $A_i$'s are disjoint) and similarly
$\sum_{j = 1}^s M_{i,j} = 0$.
Thus, the rank of $M(\A,\B)$ is less or equal to 
$\min\{s,t\}-1$. In other words
$Q(n,t,s,R) = 0$ if $R \geq \min\{s,t\}$.
 
Similarly, let $(\A,\B)$ be a pair of partitions of $[K] \times [m]$ where 
$\A =\{A_1,\hdots,A_t\} \in P(Km,t)$ (identifying $[Km]$ with $[K] \times [m]$)
and $\B = \{B_1,\hdots,B_s\} \in U^*(K,m,s)$.
Let $A_i - 1$ denote the sets 
whose elements are $\{(p,u-1),\, (p,u) \in A_i\}$.
In contrast to above we do {\it not} calculate in the circular sense this time.
So elements of the form $(p,0)$ may appear in $A_i -1$.
Then to such a pair $(\A,\B)$ we associate a matrix $L=L(\A,\B)$ with entries
\begin{equation}\label{def_LAB}
L_{i,j} \,=\, \sum_{(p,u) \in A_i \cap B_j} (-1)^p ~- \sum_{(p,u) \in (A_i-1) \cap B_j} (-1)^p.
\end{equation}
Finally, we define 
\begin{equation}\label{def_Qs}
Q^*(K,m,t,s,R) \,:=\, \# \left\{(\A,\B):\, \A \in P(Km,t), \B \in U^*(K,m,s), 
\operatorname{rank} L(\A,\B) \,=\, R\right\}.
\end{equation}
Later in Section \ref{sec_Qn} we will provide some more information
on these combinatorial quantities.

\subsection{The Main Theorems}

In order to formulate our first theorem let $F_n(\theta), n \in \N$,
denote the functions defined in terms of a generating
function by
\begin{equation}\label{def_Fn}
\sum_{n=1}^\infty F_n(\theta) \frac{x^n}{n!} \,=\, \exp(\theta(e^x-x-1)).
\end{equation}
Clearly, $F_n$ is connected to the associated Stirling numbers 
of the second kind $S_2(n,k)$ by (\ref{expgen_func}). We refer to Section
\ref{sec_Qn} for a list of $F_{2n}$ for $n=1,\hdots,6$.
Further, we define
\[
G_n(\theta) \,:=\, \theta^{-n} F_n(\theta).
\]
Also recall that $D= (2q+1)^d$.
Then our first theorem about exact reconstruction of sparse trigonometric
polynomials reads as follows.

\begin{Theorem}\label{thm_main0} 
Assume $f \in \Pi_q^d(M)$ with some sparsity $M \in \N$.
Let $x_1,\hdots,x_N \in [0,2\pi]^d$ be 
independent random variables
having the uniform distribution on $[0,2\pi]^d$.
Choose $n \in \N$, $\beta > 0,\kappa > 0$ 
and $K_1,\hdots, K_n \in \N$ such that
\begin{equation}\label{cond_kappa_0}
a \,:=\, \sum_{m=1}^n \beta^{n/K_m} < 1
\quad\mbox{ and }\quad 
\frac{\kappa}{1-\kappa} \,\leq\, \frac{1-a}{1+a}M^{-3/2}.
\end{equation}
Set $\theta:= N/M$.
Then with probability at least
\begin{equation}\label{Pbound0}
1-\left(D\beta^{-2n} \sum_{m=1}^n G_{2mK_m}(\theta) 
+ M \kappa^{-2} G_{2n}(\theta)\right)
\end{equation}
$f$ can be reconstructed exactly from
its sample values $f(x_1),\hdots,f(x_N)$ by solving 
the minimization problem (\ref{min_problem}). 
\end{Theorem}
 
We will illustrate the probability bound (\ref{Pbound0}) 
later in Section \ref{sec_bound} with some plots.
In particular, the probability of exact reconstruction is high if the 
``non-linear oversampling factor'' $\theta = N/M$ is large enough.
Of course, in order to obtain useful results one has to optimize
with respect to
the parameters occuring in (\ref{Pbound0}). In particular, the choice
of $n$ is crucial. It may not be chosen too small 
but also not too large depending
on $\theta$. Indeed, pursuing this strategy leads to the following 
qualitative result.

\begin{corollary}\label{cor1} There
exists an absolute constant $C > 0$ such that the following is true.
Assume $f \in \Pi_q^d(M)$ for some sparsity $M \in \N$.
Let $x_1,\hdots,x_N \in [0,2\pi]^d$ be 
independent random variables
having the uniform distribution on $[0,2\pi]^d$.
If for some $\epsilon > 0$ it holds
\[
N \,\geq\, C M (\log D + \log (\epsilon^{-1}))
\]
then with probability at least $1-\epsilon$ the trigonometric polynomial
$f$ can be recovered from its sample values $f(x_j), j=1,\hdots,N$,
by solving the $\ell^1$-minimization problem (\ref{min_problem}).
\end{corollary}

This formulation is similar to the main theorem in \cite{CRT1} concerned
with exact reconstruction in the context of the discrete Fourier transform.
Indeed, setting $\epsilon = D^{-\sigma}$ yields a probability of exact
reconstruction of at least $1-D^{-\sigma}$ provided 
$N \geq CM (\sigma+1) \log D$.

We remark that
(\ref{Pbound0}) of Theorem \ref{thm_main0} allows to actually compute 
precise bounds on the probability of exact reconstruction when
the parameters $M,N,D$ are given explicitly. But clearly, 
the previous corollary is easier to interpret. 
This is the reason why we have given both results.

For our next theorem we model also the set $T \subset [-q,q]^d$ 
of non-vanishing Fourier coefficients as random. So we will not treat 
arbitrary sparse polynomials, but only ``generic'' ones. The hope
is, of course, that this provides even better estimates for
the probability of exact reconstruction.

Let $0<\tau < 1$. 
The probability that an 
index $k \in [-q,q]^d$ belongs to $T$ is modelled as
\begin{equation}\label{T_model}
\P(k \in T) \,=\, \tau
\end{equation}
independently for each $k$. We also assume that the choice of $T$ 
and the choice
of the sampling set $X$ are stochastically independent.
Clearly, the expected size of $T$ is 
$\E|T| = \tau D = \tau (2q+1)^d$ and $|T|$ follows the binomial distribution.
For convenience we also introduce $\Pi_T$ as the set of all trigonometric 
polynomials whose coefficients are supported on $T$.

We also need some auxiliary notation. For $n \in  \N$ we define
\begin{equation}\label{def_W}
W(n,N,\E|T|,D) \,:=\, N^{-2n} \sum_{t=1}^{\min\{n,N\}} \frac{N!}{(N-t)!} 
\sum_{s=2}^{2n} (\E|T|)^s \sum_{R=0}^{\min\{t,s\}-1} Q(2n,t,s,R) D^{-R}
\end{equation}
and for $K,m \in \N$
\[
Z(K,m,N,\E|T|,D) \,:=\,
N^{-2Km} \sum_{t=1}^{\min\{Km,N\}} \frac{N!}{(N-t)!} 
\sum_{s=1}^{2Km} (\E|T|)^s \sum_{R=0}^{\min\{t,s\}} Q^*(2K,m,t,s,R) D^{-R}.
\]

Our second theorem about reconstructing a sparse trigonometric polynomial
from random samples by Basis Pursuit is given as follows.
 
\begin{Theorem}\label{thm_main}
Let $x_1,\hdots,x_N \in [0,2\pi]^d$ be independent random variables
having the uniform distribution on $[0,2\pi]^d$. Further
assume that $T$ is a random subset of $[-q,q]^d$ modelled
by (\ref{T_model}) (with $T$ being independent of $x_1,\hdots,x_N$)
such that $\E|T| =\tau D \geq 1$.
Choose $n \in \N$, $\alpha > 0, \beta > 0,\kappa > 0$ 
and $K_1,\hdots, K_n \in \N$ such that
\begin{equation}\label{cond_kappa_thm}
a \,:=\, \sum_{m=1}^n \beta^{n/K_m} < 1
\quad\mbox{and}\quad 
\frac{\kappa}{1-\kappa} \,\leq\, 
\frac{1-a}{1+a}\left((\alpha+1)\E|T|\right)^{-3/2}.
\end{equation}
Then with probability at least 
\begin{equation}\label{Pbound}
1 - \left(\kappa^{-2}W(n,N,\E|T|,D) 
+ \beta^{-2n} D \sum_{m=1}^n Z(K_m,m,N,\E|T|,D) 
+ \exp\left(-\frac{3\alpha^2}{6+2\alpha} \E|T|\right)\right)
\end{equation}
any $f \in \Pi_T \subset \Pi_q^d(|T|)$ can be 
reconstructed exactly from
its sample values $f(x_1),\hdots,f(x_N)$ by solving 
the minimization problem (\ref{min_problem}). 
\end{Theorem}

Of course, the theorem has to be understood in the
sense that the set $T$ is not known a priori because with
the knowledge of $T$ it would be in fact much easier to reconstruct
$f$. (Although it seems that in higher dimensions $d\geq 2$ 
not many theoretical results are available, see e.g. \cite{BG}.)

Like the previous theorem this result shows that 
the probability for reconstructing
the original sparse polynomial is indeed high for appropriate choices
for the number of sampling points $N$ and the expected sparsity $\E|T|$.
We will illustrate this later in Section \ref{sec_bound} by computing 
numerical plots for the 
bound in (\ref{Pbound}). Since the theorem does not treat arbitrary
$f$'s but only ``generic'' ones in the sense that the set $T$ is
random one may expect that the bound for the probability for
exact reconstruction is better than the one in Theorem \ref{thm_main0}.
As we will see later, this is indeed the case if one takes the same $n$,
see also Section \ref{sec_remarks}.
Unfortunately, we were not able to compute the bound explicitly
for $n\geq 5$ so that practically up to now Theorem \ref{thm_main0}
gives the better bounds since here we are to evaluate the bound 
(\ref{Pbound0}) for any $n$. 

The reason for not being able to compute (\ref{Pbound}) for arbitrary $n$
is due to the fact that we do not have an explicit expression
(or a recursion formula, or a good estimation) of the 
numbers $Q(2n,t,s,R)$ and $Q^*(2K,m,t,s,R)$.
We were only able to compute them on a computer up 
to $n=4$ by checking the rank of $M(\A,\B)$ and $L(\A,\B)$ for 
all possible pairs $(\A,\B)$ of partitions. 
Already for $n=5$ 
the computing times would exceed several days and with $n=7$ at the latest
the task nearly becomes an impossibility since the rank of 
$576 535 660 478 649 \approx 5.7\times10^{14}$ 
matrices would have to checked for computing 
the numbers $Q(14,t,s,R)$. 
So we have to leave it as
an interesting open problem to provide more information on the numbers
$Q(2n,t,s,R)$ and $Q^*(2K,m,t,s,R)$, see also Section \ref{sec_Qn}.
We hope that with a progress on this combinatorial problem
we can improve significantly our probability bounds.

\begin{remark}\label{rem_thm} 
\begin{itemize}
\item[(a)] In both theorems it is reasonable 
to choose $K_m \approx \frac{m}{n}$, for instance rounding
$m/n$ to the nearest integer. 
In this
way $m K_m \approx n$ for all $m$ and further
\[
\sum_{m=1}^n \beta^{n/K_m} \approx \sum_{m=1}^n \beta^m \,\approx\, 
\frac{\beta}{1-\beta}.  
\]
Indeed, in the limit $n \to \infty$ all the above expressions become equal.
As we require the left hand side to be less than $1$, we should choose 
$\beta$ approximately less than $1/2$. Actually, a choice near $1/2$ turned
out to be good.
\item[(b)] There is nothing special about the underlying set 
$[-q,q]^d \cap \Z^d$. Indeed, both theorem still hold when taking
any other finite subset of $\Z^d$ of size $D$ instead.
\item[(c)]
If one is interested in choosing the dimension $D=(2q+1)^d$ of
the problem very large then one may observe that
\[
\lim_{D \to \infty} W(n,N,\E|T|,D) \,=\, 
N^{-2n} \sum_{t=1}^{\min\{n,N\}} 
\frac{N!}{(N-t)!} \sum_{s=1}^{2n} Q(2n,t,s,0)\,(\E|T|)^s
\]
and
\[
\lim_{D \to \infty} \,Z(K,m,N,\E|T|,D) \,=\, N^{-2Km}
\sum_{t=1}^{\min\{Km,N\}}\frac{N!}{(N-t)!}\sum_{s=1}^{2Km} Q(2K,m,t,s,0)\,(\E|T|)^s.
\] 
(Of course, we keep $\E|T|$ fixed in this limit so that $\tau = D/\E|T|$, see
(\ref{T_model}), has to be adjusted in the process of passing with $D$ to
infinity.) This shows that the numbers $Q(2n,t,s,0)$ and 
$Q^*(2K,m,t,s,0)$ play the most important role in the probibility bound 
(\ref{Pbound}) of Theorem \ref{thm_main}. In fact, the tables in the Appendix
and Lemma \ref{lem_Qn0} indicate that these numbers are quite small for $R=0$
compared to other values of $R$.
\item[(d)] In practice, we usually do not
have precisely sparse signals. However, signals that can be approximated by
sparse ones may appear quite frequently (e.g. in the context of best $n$-term
approximation). We leave the 
investigation of related questions to future contributions, see also 
\cite{CT} for the setting of the discrete Fourier transform.
\end{itemize}
\end{remark}

\section{Proof of the Main Results}
\label{sec_Proof}

We will develop the proofs of both theorems in parallel. 
The basic idea is similar as in the paper \cite{CRT1} by Candes, Romberg
and Tao. However, there are also significant differences and, in particular,
it turns out that our approach leads to a simpler and
slightly less technical proof (although still considerably elaborate). Also
the idea of modelling the ``sparsity set'' as random is new and
requires special treatment.

Let us first introduce some auxiliary notation. By 
$\ell^2([-q,q]^d), \ell^2(T), \ell^2(X)$ we denote the $\ell^2$ space of
sequences indexed by $[-q,q]^d$, $T \subset [-q,q]^d$ and $X$, respectively, 
endowed with the usual Euclidean norm.
Moreover, we introduce the operator
\[
\F_X : \ell^2([-q,q]^d) \to \ell^2(X),
\qquad \F_X c (x_j) \,:=\, \sum_{k \in [-q,q]^d} c_k e^{i k\cdot x_j},
\quad j=1,\hdots,N.
\]
By $\F_{TX} : \ell^2(T) \to \ell^2(X)$ we denote the restriction
of $\F_X$ to sequences supported only on $T$. The
adjoint operators are denoted by $\F^*_X: \ell^2(X) \to \ell^2([-q,q]^d)$ and
$\F^*_{TX}: \ell^2(X) \to \ell^2(T)$.
 
Clearly, our problem is equivalent to reconstructing
a sequence $c$ from $\F_X c$ by solving the problem
\begin{equation}\label{min_problem2}
\min \|c'\|_1 \quad \mbox{ subject to } \quad \F_X c' \,=\, \F_X c.
\end{equation}
For $c \in \ell^2([-q,q]^d)$ we introduce its sign by 
\[
\sgn(c)_k \,=\, \frac{c_k}{|c_k|},  \quad k \in \supp c,
\quad\mbox{and}\quad \sgn (c)_k \,=\, 0,\quad k\notin \supp c.
\]
Hereby, $\supp c$ denotes the support of $c$.

The key lemma for our proofs is the following
duality principle.

\begin{lemma}\label{lem_key}
Let $c \in \ell^2([-q,q]^d)$ and $T:=\supp c$.  
Assume $\F_{TX}: \ell^2(T) \to \ell^2(X)$ to be injective. 
Suppose that there
exists a vector $P \in \ell^2([-q,q]^d)$ with the following properties:
\begin{itemize}
\item[(i)] $P_k = \sgn c_k$ for all $k \in T$,
\item[(ii)] $|P_k| < 1$ for all $k \notin T$,
\item[(iii)] there exists a vector $\lambda \in \ell^2(X)$ such that
$P= \F_X^* \lambda$.
\end{itemize}
Then $c$ is the unique minimizer to the problem 
(\ref{min_problem2}).
\end{lemma}
\begin{Proof} The proof mimiques the one by Candes, Romberg and 
Tao \cite[Lemma 2.1]{CRT1}. For the sake of completeness we repeat the
argument.

Let $b$ be a vector with $\F_X b = \F_X c$ and set $h:= b-c$. Clearly,
$\F_X h$ vanishes. For any $k \in  T$ we have
\[
|b_k| \,=\, |c_k + h_k| \,=\, ||c_k| + h_k \ol{\sgn(c)_k}|
\,\geq \, |c_k| + \Re(h_k \ol{\sgn(c)_k})
\,=\, |c_k| + \Re(h_k \ol{P_k}).
\]  
If $k \notin T$ then $|b_k| = |h_k| \geq \Re(h_k \ol{P_k})$ since
$|P_k| < 1$. This gives
\[
\|b\|_1 \,\geq\, \|c\|_1 + \sum_{k \in [-q,q]^d} \Re(h_k \ol{P_k}).
\]
Further, observe that
\[
\sum_{k \in [-q,q]^d} \Re(h_k \ol{P_k})
\,=\, \Re\left( \sum_{k \in [-q,q]^d} h_k \ol{(\F^*_X \lambda)_k}\right)
\,=\, \Re \left( \sum_{j = 1}^N (\F_X h)(x_j) \ol{\lambda(x_j)} \right)
\,=\, 0
\]
since $\F_X h$ vanishes. Altogether, we proved $\|b\|_1 \geq \|c\|_1$,
and thus $c$ is a minimizer of (\ref{min_problem2}). 

It remains to prove the uniqueness. The above argument shows
that having the equality $\|b\|_1 = \|c\|_1$ forces 
$|h_k| = \Re(h_k \ol{P_k})$ for all $k \notin T$. Since $|P_k|<1$
this means that $h$ vanishes outside $T$. Since also 
$\F_X h$ vanishes, it follows from
the injectivity of $\F_{TX}$ that $h$ vanishes identically and hence,
$b=c$. This shows that $c$ is the unique minimizer of (\ref{min_problem2}).
\end{Proof}

Concerning the assumption on the injectivity of $\F_{TX}$ we have
the following simple result. 

\begin{lemma}\label{lem_injective} If $N \geq |T|$ 
then $\F_{TX}$ is injective almost
surely. 
\end{lemma}
\begin{Proof} The proof is essentially contained in \cite[Theorem 3.2]{BG}. 
There it is proved that any $|T| \times |T|$ 
submatrix of $\F_{TX}$ has non-vanishing
determinant almost surely (even under slightly more general assumptions
on the distribution of the random variables $x_1,\hdots,x_N$).
This implies the result.
\end{Proof} 

Now our strategy for proving Theorem
\ref{thm_main} is obvious. We need 
to show that with high probability there exists a vector $P$ with
the properties assumed in Lemma \ref{lem_key}.
To this end we proceed similarly as in \cite{CRT1}. (Actually, the injectivity of $\F_{TX}$
will also follow from this finer analysis so that Lemma \ref{lem_injective} will not
be needed in the end.)

We introduce the restriction operator 
$R_T : \ell^2([-q,q]^d) \to \ell^2(T)$, $R_T c_k = c_k$ for $k \in T$.
Its adjoint $R_T^*= E_T : \ell^2(T) \to \ell^2([-q,q]^d)$ 
is the operator that extends a vector outside
$T$ by zero, i.e., $(E_T d)_k = d_k$ for $k \in T$ 
and $(E_T d)_k = 0$ otherwise.

Now assume for the moment that 
$\F_{TX}^* \F_{TX} : \ell^2(T) \to \ell^2(T)$ is
invertible. 
(By Lemma \ref{lem_injective} this is true almost surely if $N \geq |T|$
since $\F_{TX}$ is then injective.) 
In this case we define $P$ explicitly by
\[
P\,:=\, \F^*_X \F_{TX} (\F_{TX}^* \F_{TX})^{-1} R_T \sgn(c),
\]
where as before $T= \supp c$. Then clearly $P$ has property (i)
and property (iii) in Lemma \ref{lem_key} with 
\[
\lambda:= \F_{TX} (\F_{TX}^* \F_{TX})^{-1} R_T \sgn(c) \in \ell^2(X).
\]
We are left with proving that $P$ has property (ii) of Lemma \ref{lem_key}
with high probability.

To this end we introduce the auxiliary operators
\[
H : \ell^2(T) \to \ell^2([-q,q]^d), \qquad
H \,:=\, N E_T - \F^*_X \F_{TX}
\]
and
\[
H_0 : \ell^2(T) \to \ell^2(T),\qquad
H_0 \,:=\, R_T H \,=\, N I_T - \F^*_{TX} \F_{TX} ,
\]
where $I_T$ denotes the identity on $\ell^2(T)$. Obviously, $H_0$
is self-adjoint, and $H$ acts on
a vector as
\[
(H c)_{\ell} \,=\, - \sum_{j=1}^N \sum_{\substack{k\in T \\ k\neq \ell}} c_k e^{i(k-\ell)\cdot x_j}.
\]
Now we can write
\[
P \,=\, (NE_T - H) \left(N I_T -H_0\right)^{-1} R_T \sgn(c).
\]
As we are interested in property (ii) in Lemma \ref{lem_key} we consider
only values of $P$ on $T^c = [-q,q]^d \setminus T$. Since $R_{T^c} E_T = 0$ 
we have
\[
P_k \,=\, -\frac{1}{N} R_{T^c} H (I_T - \frac{1}{N} H_0)^{-1} R_T \sgn(c) 
\quad \mbox{ for all }
k \in T^c.
\]
Let us look closer at the term $(I_T - \frac{1}{N} H_0)^{-1}$. 
To this end let $n \in \N$ be some arbitrary number.
By the von Neumann series we can write
\[
\left(I_T - \left(\frac{1}{N} H_0\right)^n\right)^{-1} \,=\, I_T + A_n
\]
with
\begin{equation}\label{def_An}
A_n \,:=\, \sum_{r=1}^\infty (\frac{1}{N} H_0)^{rn}.
\end{equation}
Using the identity 
\begin{equation}\label{pow_inv}
(1-M)^{-1} \,=\, (1- M^n)^{-1}(1 + M + \cdots + M^{n-1})
\end{equation}
we obtain
\[
(I_T - \frac{1}{N}H_0)^{-1} \,=\, (I_T +A_n)\sum_{m=0}^{n-1} 
\left(\frac{1}{N} H_0\right)^m.
\]
Thus, on the complement of $T$, we may write
\begin{eqnarray*}
R_{T^c} P &=&  -\frac{1}{N} H (I_T + A_n) 
\left(\sum_{m=0}^{n-1} (N^{-1}H_0)^m\right) R_T \sgn(c)\\
&=& - \sum_{m=1}^n (N^{-1} H R_T)^m \sgn(c)
-\frac{1}{N} H A_n R_T \sum_{m=0}^{n-1} (N^{-1} HR_T)^m \sgn(c)\\
&=& -(P^{(1)} + P^{(2)}),
\end{eqnarray*}
where
\[
P^{(1)} \,=\, S_n \sgn(c), \quad\mbox{ and }\quad 
P^{(2)} \,=\, \frac{1}{N} HA_n R_T (I + S_{n-1}) \sgn(c),
\]
with
\[
S_n \,:=\, \sum_{m=1}^n (N^{-1} H R_T)^m.
\]
Our aim is to estimate $\P(\sup_{k \in T^c} |P_k| \geq 1)$. To this end let
$a_1,a_2>0$ be numbers satisfying $a_1 + a_2 = 1$. Then
\begin{equation}\label{Psplit}
\P(\sup_{k \in T^c} |P_k| \geq 1) \,\leq\, 
\P\left(\{\sup_{k \in T^c} |P^{(1)}_k| \geq a_1\} \cup 
\{\sup_{k \in T^c} |P^{(2)}_k| \geq a_2\}\right).
\end{equation}
Clearly,
\begin{align}
\P(|P^{(1)}_k|\geq a_1) \,& = \, 
\P\left(\left|\left(\sum_{m=1}^n (N^{-1} HR_T)^m \sgn(c)\right)_k\right|\geq a_1\right)\notag\\
\label{Pinf_estim}
& \leq \, \P\left(\sum_{m=1}^n |((N^{-1}HR_T)^m \sgn(c))_k| \geq a_1\right)
\,=:\, \P(E_k).
\end{align}
Consider $P^{(2)}$.  Denoting
$\ell^\infty = \ell^\infty([-q,q]^d)$ the space of sequences indexed by
$[-q,q]^d$ with the supremum norm (and similarly defining $\ell^\infty(T)$) 
we have
\begin{equation}\label{P2_ineq}
\sup_{k \in T^c} |P^{(2)}_k| \leq \|P^{(2)}\|_\infty
\leq \|\frac{1}{N}H A_n\|_{\ell^\infty(T) \to \ell^\infty} 
(1 + \|R_T S_{n-1} \sgn(c)\|_{\ell^\infty(T)})
\end{equation}
In order to analyze the term $\|R_T S_{n-1} \sgn(c)\|_{\ell^\infty(T)}$
we observe that similarly as in (\ref{Pinf_estim})
\[
\P(|(S_{n-1} \sgn(c))_k| \geq a_1)\,\leq\, 
\P\left(\sum_{m=1}^n |((N^{-1}HR_T)^m \sgn(c))_k| \geq a_1\right) \,=\,\P(E_k).
\]
Let us now treat the operator norm appearing in (\ref{P2_ineq}).
For simplicity we write $\|\cdot\|_\infty$ instead of 
$\|\cdot\|_{\ell^\infty \to \ell^\infty}$. It holds 
$\|A\|_\infty = \sup_r \sum_s |A_{rs}|$. Clearly, 
\[
\|\frac{1}{N} H A_n\|_{\ell^\infty} \leq \|\frac{1}{N} H\|_\infty 
\|A_n\|_{\ell^\infty(T)}.
\]
Moreover, $\|\frac{1}{N}H\|_\infty \leq |T|$ as $H$ has $|T|$ columns and 
each entry is bounded by $N$ in absolute value.

In order to analyze $A_n$ we will work with the Frobenius norm. 
For a matrix $A$ it is defined as
\[
\|A\|_F^2 \,:=\, \Tr(A A^*) \,=\, \sum_{r,s} |A_{rs}|^2
\]
where $\Tr(A A^*)$ denotes the trace of $A A^*$. 
Assume for the moment that 
\begin{equation}\label{assume_H0}
\|(\frac{1}{N} H_0)^n\|_F \,\leq\, \kappa < 1.
\end{equation}
Then it follows directly from the definition (\ref{def_An}) of  $A_n$ that
\[
\|A_n\|_F \,=\, \left\|\sum_{r=1}^\infty \left(\frac{1}{N} H_0\right)^{rn}\right\|_F
\,\leq\, \sum_{r=1}^\infty \| (N^{-1} H_0)^n\|_F^r 
\,\leq\, \sum_{r=1}^\infty \kappa^{r} 
\,=\, \frac{\kappa}{1-\kappa}.
\]
Moreover, since $A_n$ has $|T|$ columns 
it follows from the Cauchy-Schwarz inequality that
\[
\|A_n\|_\infty^2 \,\leq\, \sup_{i} |T| \sum_{j} |A_n(i,j)|^2
\,\leq\, |T| \|A_n\|_F^2.
\] 
So assuming (\ref{assume_H0})
and $\|S_{n-1} \sgn(c)\|_\infty < a_1$ we have
\[
\sup_{t \in T^c} |P^{(2)}_k| \,\leq\, (1+a_1) |T|^{3/2} 
\frac{\kappa}{1-\kappa}.
\]
In particular, if 
\begin{equation}\label{cond_kappa00}
\frac{\kappa}{1-\kappa} \leq \frac{a_2}{1+a_1}|T|^{-3/2}
\end{equation}
then $ \sup_{t\in T^c} |P^{(2)}_k| \leq a_2$ as desired.
Also it follows from (\ref{cond_kappa00}) that $\kappa < 1$ as
$|T| \geq 1$ without loss of generality (if $T=\emptyset$ then $f=0$
and $\ell^1$-minimization will clearly recover $f$.)

Now we have to distinguish between the situation in Theorem 
\ref{thm_main0} and the one in Theorem \ref{thm_main} since 
in the latter $|T|$ is a random variable while in the first 
it is deterministic.

\begin{enumerate}
\item 
Let us first treat the case of Theorem \ref{thm_main} where $|T|$ is random.
If 
\[
|T| \,\leq\, (\alpha+1) \E|T|
\]
with $\alpha > 0$ and
\begin{equation}\label{cond_kappa}
\frac{\kappa}{1-\kappa} \leq \frac{a_2}{1+a_1} ((\alpha+1)\E|T|)^{-3/2}
\end{equation}
then clearly (\ref{cond_kappa00}) is satisfied and consequently
\[
\sup_{t \in T^c} |P_k^{(2)}| \leq a_2.
\]
Using the union bound we altogether obtain from (\ref{Psplit})
\begin{align}
&\P(\sup_{k \in T^c} |P_k| \geq 1) 
\leq
\P\left(\bigcup_{k \in T^c}\{|P_k^{(1)}| \geq a_1\} 
\cup \{\|R_T \sgn(c)\|_{\ell^\infty(T)} \geq a_1\} 
\cup \{\|(N^{-1}H_0)^n\|_F \geq \kappa\}\right. \notag\\
&\left.\phantom{\P(\sup_{k \in T^c} |P_k| \geq 1) \leq \P ( \bigcup_{k \in T^c})}
\cup \{|T| \geq (\alpha + 1) \E|T| \}\right).\notag\\
&\leq \P\left(\bigcup_{k \in [-q,q]^d} E_k \cup  
\{\|(N^{-1}H_0)^n\|_F \geq \kappa\} 
\cup \{|T| \geq (\alpha + 1) \E|T|\}\right)\notag\\
\label{Pestim_3}
&\leq \sum_{k \in [-q,q]^d} \P(E_k) + 
\P(\|(N^{-1}H_0)^n\|_F \geq \kappa)
+ \P(|T| \geq (\alpha + 1) \E|T|).
\end{align}
As $|T|$ is the sum of independent random variables we obtain for the third
term from the large deviation theorem (see for instance 
equation (6) in \cite{Boucheron}, where
also slightly better estimates are available) 
\begin{align}
\P(|T| \geq \E|T| + \alpha \E|T|) \,&\leq\, 
\exp\left(-(\alpha \E|T|)^2/(2\E|T| + 2(\alpha \E|T|)/3)\right) \notag\\
\label{prob_ET}
&=\, \exp\left(-\frac{3 \alpha^2}{6+2\alpha}\E|T|\right).
\end{align}
So we are left with the two other expressions in (\ref{Pestim_3}).

\item In the situation of Theorem \ref{thm_main0} we proceed in almost the same
way with the only difference that we do not need to treat $|T|$ as random
variable. Under the condition in (\ref{cond_kappa00}) this yields
\begin{equation}\label{Pestim30}
\P(\sup_{k \in T^c} |P_k| \geq 1) \,\leq\, \sum_{k \in [-q,q]^d} \P(E_k) + 
\P(\|(N^{-1}H_0)^n\|_F \geq \kappa).
\end{equation}
Hence, also here we need to estime $\P(E_k)$ and 
$\P (\|(N^{-1}H_0)^n\|_F \geq \kappa)$.
\end{enumerate}

\subsection{Analysis of powers of $H_0$}
\label{sec_H0}

In this section we treat the second term in (\ref{Pestim_3}) 
and (\ref{Pestim30}), i.e.,
we estimate powers of the random matrix $H_0$ in the Frobenius norm.
To this end Markov's inequality suggests to estimate the expectation
of $\|H_0^n\|^2_F$. In the following lemma we only take the expectation
with respect to the random sampling set $X=\{x_1,\hdots,x_N\}$.
For the situation of Theorem \ref{thm_main} we postpone the computation
of the full expectation $\E = \E_T \E_X$ (the latter by Fubini's theorem).

\begin{lemma}\label{lem_H0} It holds
\[
\E_X\left[\|H_0^{n}\|^2_F\right] \,=\, 
\sum_{t=1}^{\min\{n,N\}} \frac{N!}{(N-t)!}
\sum_{\A \in P(2n,t)} \sum_{\substack{k_1,\hdots,k_{2n} \in T\\ k_j \neq k_{j+1}, j \in [2n]}}
\prod_{A\in \A} \delta\left(\sum_{r\in A} (k_{r+1}-k_r)\right) 
\]
where $\delta(n)$ denotes the Kronecker $\delta_{0n}$ and 
$k_{2n+1} = k_1$.
\end{lemma} 
\begin{Proof}
As $H_0$ is self-adjoint we need to estimate $\|H_0^n\|_F^2 = \Tr(H_0^{2n})$.
Observe that 
\[
H_0(k,k') \,=\, h(k'-k), \qquad k,k' \in T
\]
with
\[
h(k) \,=\, -\delta(k) \sum_{\ell=1}^N e^{ik\cdot x_\ell}.
\]
Thus, $H_0^2(k,k') \,=\, \sum_{t \in T, t\neq k,k'} h(t-k)h(k'-t)$
and
\[
H_0^{2n}(k_1,k_1) \,=\, \sum_{k_2,\hdots,k_{2n} \in T, k_j \neq k_{j+1}} h(k_2-k_1)
\cdots h(k_{1}-k_{2n})
\]
where we agree on the convention that $k_{2n+1} = k_1$. This yields
\[
\Tr(H_0^{2n}) \,=\, \sum_{k_1,\hdots,k_{2n} \in T, k_j \neq k_{j+1}} h(k_2-k_1) h(k_3-k_2)
\cdots h(k_1-k_{2n}).
\]
Using linearity of expectation and the definition of $h$ we get
\[
\E_X\left[\Tr(H_0^{2n})\right] 
\,=\, \sum_{\ell_1,\hdots,\ell_{2n} = 1}^N ~
\sum_{\substack{k_1,\hdots,k_{2n} \in T\\ k_j \neq k_{j+1}}}
\E_X\left[ \exp\left(i \sum_{r=1}^{2n} (k_{r+1} - k_{r})\cdot x_{\ell_r}\right)\right].
\]
Let us consider the latter expected value. Here we have
to take into accound that some of the indeces $\ell_r$ might be the same.
This is where set partitions enter the game.

We associate a partition $\A = (A_1,\hdots,A_t)$ of $\{1,\hdots,2n\}$
to a certain vector $(\ell_1,\hdots,\ell_{2n})$ such that
$\ell_r = \ell_{r'}$ if and only if $r$ and $r'$ are contained in the same
set $A_i \in \A$. This is allows us to unambiguously write $\ell_A$ 
instead of $\ell_r$ if $r \in  A$. The independence of the 
$x_{\ell_A}$ yields
\begin{align}\label{ExpectX}
\E_X & \left[\exp\left(i \sum_{r=1}^{2n} (k_{r+1}-k_r)\cdot x_{\ell_r} \right)\right]
\,=\, \E_X\left[\exp\left(i \sum_{A \in \A} \sum_{r \in A} (k_{r+1}-k_r) \cdot x_{\ell_A}\right)\right]\notag\\
&=\, \prod_{A \in \A} \E_X\left[\exp\left(i \sum_{r\in A} (k_{r+1}-k_r)\cdot x_{\ell_A}
\right)\right].
\end{align}
Since $x_{\ell_A}$ has the uniform distribution on $[0,2\pi]^d$ we
obtain
\begin{align}
\E_X \left[ \exp\left(i \sum_{r \in A}(k_{r+1} - k_{r})\cdot x_{\ell_A}\right) \right]
\, & = \, \int_{[0,2\pi]^d} \exp\left(i  \sum_{r \in A} (k_{r+1} - k_r) \cdot x\right) dx\notag\\
\label{Expect_Trig}
& = \, \delta\left( \sum_{r \in A} (k_{r+1} - k_r)\right).
\end{align}
Observe that the last expression 
is independent of 
the precise values of the $\ell_r$. Only the generated partition $\A$ plays a role. 
Moreover, if $A \in \A$ contains only one element then 
(\ref{Expect_Trig}) vanishes
due to the condition $k_{r+1} \neq k_r$. Thus, we only
need to consider partitions $\A$ satisfying $|A|\geq 2$ for 
all $A \in \A$, i.e., partitions in $P(2n,t)$.
Moreover, observe 
that the number of vectors $(\ell_{A_1}, \hdots, \ell_{A_t}) \in \{1,\hdots,N\}^t$
with different entries is 
precisely $N\cdots(N-t+1) = N! / (N-t)!$ if $N\geq t$ and $0$ if 
$N\leq t$.
Finally, we obtain
\[
\E_X[\|H_0^n\|_F^2] \,=\, \sum_{t=1}^{\min\{n,N\}} \sum_{\A \in P(2n,t)}
\sum_{\substack{k_1,\hdots,k_{2n} \in T\\ k_j \neq k_{j+1}}}
\prod_{A\in \A} \delta\left(\sum_{r\in A} (k_{r+1}-k_r)\right),
\]
which is precisely the content of the lemma.
\end{Proof}

In view of the previous lemma we define
for simplicity and later reference 
\begin{align}\label{def_CAT}
&C(\A,T) \,:=\, \sum_{\substack{k_1,\hdots,k_{2n} \in T\\ k_j \neq k_{j+1}}}
\prod_{A\in \A} \delta\left(\sum_{r\in A} (k_{r+1}-k_r)\right)\\
&=\, \#\left\{ (k_1,\hdots,k_{2n}) \in T^{2n}:\, k_j \neq k_{j+1}, j \in [2n],
\mbox{ and }
\sum_{r \in A} (k_{r+1} - k_r) = 0 \mbox{ for all } A \in \A \right\}.\notag
\end{align}

\subsection{Analysis of $\P(E_k)$}

Let us now treat the first term $\P(E_k)$ in (\ref{Pestim_3}) resp. 
(\ref{Pestim30}). 
To this end let $\beta_m, m=1,\hdots,n$, be positive numbers satisfying 
\[
\sum_{m=1}^n \beta_m = a_1
\]
and $K_m \in \N$, $m=1,\hdots,n$,
some natural numbers. Let $k \in [-q,q]^d$.
Using Markov's inequality in the last step we obtain
\begin{align}
\P(E_k) \,& = \, \P\left(\sum_{m=1}^n |((N^{-1}HR_T)^m \sgn(c))_k| \geq a_1\right)
\,\leq\, \sum_{m=1}^n \P(N^{-m} |((HR_T)^m \sgn(c))_k| \geq \beta_m) \notag\\
&=\, \sum_{m=1}^n \P\left(N^{-2mK_m} |((HR_T)^m \sgn(c))_k|^{2K_m} \geq \beta_m^{2K_m}\right)\notag\\
\label{P1_estim}
&\leq \, \sum_{m=1}^n \E\left[|((HR_T)^m\sgn(c))_k|^{2K_m}\right] 
N^{-2mK_m} \beta_m^{-2K_m}.
\end{align}
Let us choose $\beta_m = \beta^{n/K_m}$, i.e., 
$\beta_m^{-2K_m} = \beta^{-2n}$. This yields
\begin{equation}\label{prob_Ek}
\P(E_k) \,\leq\, \beta^{-2n} \sum_{m=1}^n 
\E\left[|((HR_T)^m\sgn(c))_k|^{2K_m}\right] N^{-2mK_m}
\end{equation}
and the condition $a_1=\sum_{m=1}^n \beta_m$ reads
\[
a_1 \,=\, a \,=\, \sum_{m=1}^n \beta^{n/K_m} < 1.
\] 
The following lemma is concerned with the expectation appearing in (\ref{prob_Ek}). 
We first investigate the expectation with respect to $X$. 
The following proof is similar
to the one of Lemma \ref{lem_H0}.

\begin{lemma}\label{lem_Expect1} 
For $k \in [-q,q]^d$ and $c \in \ell^2([-q,q]^d)$ with
$\supp c = T$ we have
\begin{align}
&\E_X\left[|((HR_T)^m\sgn(c))_{k}|^{2K}\right] \notag\\
&\leq \, 
\sum_{t=1}^{\min\{Km,N\}} \frac{N!}{(N-t)!} 
\sum_{\A \in P(2Km,t)} 
\sum_{\substack{ k^{(1)}_1,\hdots,k^{(1)}_m \in T\\
      \vdots\\
      k^{(2K)}_1,\hdots, k^{(2K)}_m \in T\\
      k^{(p)}_{j-1} \neq k^{(p)}_{j},\,j\in [m] }}
\prod_{A\in \A} 
\delta\left(\sum_{(r,p) \in A} (-1)^p(k_r^{(p)} -k_{r-1}^{(p)})\right)\notag
\end{align}
with $k^{(p)}_0 := k$ for $p=1,\hdots,2K$. Hereby, we identify partitions
of $[2Km]$ in $P(2Km,t)$ with partitions of $[2K]\times[m]$ in an obvious way.
\end{lemma}
\begin{Proof} Set $\sigma := \sgn(c)$. An elementary calculation yields
\[
((HR_T)^m \sigma)_{k} \,=\, 
(-1)^m \sum_{\ell_1,\hdots,\ell_m=1}^N 
\sum_{\substack{ k_1,\hdots,k_m \in T\\
                                          k_{j-1} \neq k_{j}, j=1,\hdots,m}} 
\sigma(k_m) e^{i(k_m-k_{m-1})\cdot x_{\ell_m}} \cdots e^{i(k_1-k_0)\cdot x_{\ell_1}}
\]
with $k_0 := k$. Thus,
\begin{eqnarray*}
|((HR_T)^m \sigma)_{k_0}|^{2} &=& 
 \sum_{\ell^{(1)}_1,\hdots,\ell^{(1)}_m=1}^N
\sum_{\ell^{(2)}_1, \hdots, \ell^{(2)}_m=1}^N
\sum_{\substack{ {k^{(1)}_1},\hdots,{k^{(1)}_m} \in T\\
     {k^{(2)}_1},\hdots, {k^{(2)}_m} \in T\\
     {k^{(p)}_{j-1}} \neq {k^{(p)}_{j}}, j\in [m], p=1,2}}   
\sigma(k_m^{(1)}) \ol{\sigma(k_m^{(2)})} \times \\
&\times & e^{i\sum_{r=1}^m(k^{(1)}_r-k^{(1)}_{r-1})\cdot x_{\ell^{(1)}_r}}
e^{-i\sum_{r=1}^m  (k^{(2)}_r-k^{(2)}_{r-1})\cdot x_{\ell^{(2)}_r}}
\end{eqnarray*}
where $k^{(1)}_0 = k^{(2)}_0 = k_0 = k$. Taking a $2K$-th power yields
\begin{align}
|((HR_T)^m \sgn(c))_{k}|^{2K}\,
&=\, \sum_{\substack{ \ell^{(1)}_1,\hdots,\ell^{(1)}_m = 1\\
           \vdots\\
          \ell^{(2K)}_1,\hdots, \ell^{(2K)}_m = 1}}^N 
\sum_{\substack{ k^{(1)}_1,\hdots,k^{(1)}_m \in T\\
       \vdots\\
       k^{(2K)}_1,\hdots, k^{(2K)}_m \in T\\
       k^{(p)}_{j-1} \neq k^{(p)}_{j} }}
\sigma(k_m^{(1)}) \ol{\sigma(k_m^{(2)})} \cdots \sigma(k^{(2K-1)}_m)\ol{\sigma(k^{(2K)}_m)}\times\notag\\
&\times \, \exp\left(i\sum_{p=1}^{2K} (-1)^p \sum_{r=1}^m (k_r^{(p)} - k_{r-1}^{(p)})\cdot x_{\ell^{(p)}_r}\right)\notag
\end{align}
with $k_0^{(p)} = k$, $p=1,\hdots,2K$.
Further, recall that $|\sigma(k)|=1$ on $T$. Taking the expected
value $\E_X$ yields
\begin{multline}
\E_X\left[|((HR_T)^m \sgn(c))_{k}|^{2K}\right]\\
\label{Expect1}
\,\leq\,
\sum_{\substack{ \ell^{(1)}_1,\hdots,\ell^{(1)}_m = 1\\
                 \vdots\\
         \ell^{(2K)}_1,\hdots, \ell^{(2K)}_m = 1}}^N 
\sum_{\substack{ k^{(1)}_1,\hdots,k^{(1)}_m \in T\\
                 \vdots\\
                k^{(2K)}_1,\hdots, k^{(2K)}_m \in T\\
                k^{(p)}_{j-1} \neq k^{(p)}_{j} }}
\E_X\left[ \exp\left(i\sum_{p=1}^{2K} (-1)^p \sum_{r=1}^m (k_r^{(p)} - k_{r-1}^{(p)})\cdot x_{\ell^{(p)}_r}\right)\right]
\end{multline}
(with equality if all the entries of $\sigma$ are equal on $T$).

Let us consider the expected value appearing in the sum.
As in the proof of Lemma \ref{lem_H0} we have to take into account that some of the indeces
$\ell^{(p)}_r$ might coincide. This affords to introduce some additional
notation. Let $(\ell^{(p)}_r)_{r=1,\hdots,m}^{p=1,\hdots,2K} \subset \{1,\hdots,N\}^{2Km}$ 
be some vector of indeces and let 
$\A = (A_1,\hdots, A_t)$, $A_i \subset \{1,\hdots,m\} \times \{1,\hdots 2K\}$ 
be a corresponding partition 
such that $(r,p)$ and $(r',p')$ are contained
in the same block if and only if  $\ell^{(p)}_{r} = \ell^{(p')}_{r'}$. For some
$A \in \A$ we may unambigously write $\ell_A$ instead of $\ell_r^{(p)}$ if $(r,p) \in A$.

Like in (\ref{ExpectX}), using that all $\ell_A$ for $A\in \A$ 
are different and that the $x_{\ell_A}$ are
independent we may write the expectation in the sum in (\ref{Expect1}) as
\begin{align}
\E & \left[  \exp\left(i\sum_{p=1}^{2K} (-1)^p \sum_{r=1}^m (k_r^{(p)} - k_{r-1}^{(p)})\cdot x_{\ell^{(p)}_r}\right)\right]\notag\\
& = \, \prod_{A \in \A} \E\left[\exp\left(i \sum_{(r,p) \in A} (-1)^p(k_r^{(p)} - k_{r-1}^{(p)})\cdot x_{\ell_A}\right)\right] 
\,=\, \prod_{A \in \A} \delta\left( \sum_{(r,p) \in A} (-1)^p(k_r^{(p)} - k_{r-1}^{(p)})\right).\notag
\end{align}
Once again,
if $A \in \A$ contains only one element then the last expression 
vanishes due to the condition 
$k_r^{(p)} \neq k_{r-1}^{(p)}$. Thus, we only
need to consider partitions $\A$ 
in $P(2Km,t)$.
Now we are able to rewrite the inequality in (\ref{Expect1}) as
\begin{align}
\E_X&\left[|((HR_T)^m \sigma)_{k}|^{2K}\right]\notag\\
\,&\leq\, \sum_{t=1}^{Km} \sum_{\A \in P(2Km,t)} 
\sum_{\substack{\ell_{(1)},\hdots,\ell_{(t)}=1 \\ \ell_{(1)},\hdots,\ell_{(t)} \mbox{\small p.w. different}}}^N
 \sum_{\substack{ k^{(1)}_1,\hdots,k^{(1)}_m \in T\\
            \vdots\\
         k^{(2K)}_1,\hdots, k^{(2K)}_m \in T\\
         k^{(p)}_{j-1} \neq k^{(p)}_{j} }}
\prod_{A\in \A} \delta\left(\sum_{(r,p) \in A} (-1)^p(k_r^{(p)} -k_{r-1}^{(p)})\right)
\notag\\
&=\, \sum_{t=1}^{\min\{Km,N\}} \frac{N!}{(N-t)!} \sum_{\A \in P(2Km,t)}
 \sum_{\substack{ k^{(1)}_1,\hdots,k^{(1)}_m \in T\\
                   \vdots\\
                 k^{(2K)}_1,\hdots, k^{(2K)}_m \in T\\
                 k^{(p)}_{j-1} \neq k^{(p)}_{j} }}
\prod_{A\in \A} 
\delta\left(\sum_{(r,p) \in A} (-1)^p(k_r^{(p)} -k_{r-1}^{(p)})\right).\notag
\end{align}
This proves the lemma.
\end{Proof}

In view of the previous lemma and
for the sake of simple notation we denote
\begin{equation}\label{def_BAT}
B(\A,T) \,:=\,  \sum_{\substack{ k^{(1)}_1,\hdots,k^{(1)}_m \in T\\
                   \vdots\\
                 k^{(2K)}_1,\hdots, k^{(2K)}_m \in T\\
                 k^{(p)}_{j-1} \neq k^{(p)}_{j}, j \in [m], p\in [2K] }}
\prod_{A\in \A} 
\delta\left(\sum_{(r,p) \in A} (-1)^p(k_r^{(p)} -k_{r-1}^{(p)})\right).
\end{equation}

\subsection{Proof of Theorem \ref{thm_main0}}

Let us assemble all the pieces to complete the proof of 
Theorem \ref{thm_main0}. By Lemma \ref{lem_H0} we need to investigate
the quantity $C(\A,T)$ defined in (\ref{def_CAT})
for $\A \in P(2n,t)$.
Here the indeces $(k_1,\hdots,k_{2n}) \in T^{2n}$ 
are subjected to the $|\A|=t$ linear
constraints $\sum_{r \in A}({k_{r+1}-k_r}) = 0$ for all $A \in \A$. These
constraints are independent except for $\sum_{r=1}^{2n} (k_{r+1}-k_r) = 0$.
Thus, we can estimate 
\begin{equation}\label{CAT_deterministic}
C(\A,T) \,\leq\, |T|^{2n-t+1} \,\leq\, M^{2n-t+1}. 
\end{equation}
By Lemma
\ref{lem_H0} we obtain (note that in the situation of Theorem \ref{thm_main0}
$T$ is not random, so $\E=\E_X$) 
\[
\E[\|H_0^n\|_F^2]
\,\leq\, \sum_{t=1}^{\min\{n,N\}} \frac{N!}{(N-t)!} 
\sum_{\A \in P(2n,t)} |T|^{2n+1-t}
\,\leq\, M^{2n+1} \sum_{t=1}^n (N/M)^t S_2(2n,t).
\]   
where $S_2(n,t)=|P(2n,t)|$ are the associated Stirling 
numbers of the second kind.
Set $\theta = N/M$.
From the generating function (\ref{expgen_func}) of the numbers
$S_2(n,k)$ we know that
\[
\sum_{t=1}^n S_2(2n,t) \theta^t \,=\, F_{2n}(\theta)
\]
with $F_{2n}$ defined by (\ref{def_Fn}). Markov's inequality
yields
\begin{align}
\P(\|(N^{-1} H_0)^n\|_F \geq \kappa)
\,&=\, \P(\|H_0^n\|_F^2 \geq N^{2n} \kappa^2)\notag\\
&\leq\, N^{-2n} \kappa^{-2} \E[\|H_0^n\|_F^2] 
\,\leq\, \kappa^{-2} M \theta^{-2n} F_{2n}(\theta) \,=\, 
\kappa^{-2} M G_{2n}(\theta). \notag
\end{align}
We remark that by (\ref{assume_H0}) we have 
$\kappa < 1$. 
In the event
that $\|(N^{-1} H_0)^n\|_F \leq \kappa$ this implies
that $(I_T - (N^{-1} H_0)^n)$ is invertible by the von Neumann series
and by (\ref{pow_inv}) also 
\[
\F^*_{TX} \F_{TX} \,=\, N(I_T - N^{-1} H_0)
\]  
is invertible. In particular, $\F_{TX}$ is injective. So this basic condition
in Lemma \ref{lem_key} is satisfied automatically  
with a probability that can be derived from the estimation above, and we
do not even need to invoke Lemma \ref{lem_injective}.

Let us now consider $\P(E_k)$. By Lemma \ref{lem_Expect1} we need to
bound $B(\A,T)$ defined in (\ref{def_BAT}), i.e., 
the number of vectors $(k_j^{(p)}) \in T^{2Km}$ satisfying
$\sum_{(r,p) \in A} (-1)^p(k_r^{(p)} -k_{r-1}^{(p)})=0$ for all $A \in \A$
with $\A \in P(2Km,t)$. These are $t$ independent linear constraints.
So the number of these indices is bounded from above by
$|T|^{2Km-t} \leq M^{2Km-t}$. Thus, similarly as above we obtain
\[
\E\left[|((HR_T)^m \sgn(c))_k|^{2K}\right]
\,\leq\,\sum_{t = 1}^{Km} N^t S_2(2Km,t) M^{2Km-t}
\,=\, M^{2Km} F_{2Km}(\theta). 
\]
By (\ref{prob_Ek}) this yields
\[
\P(E_k) \,\leq\, 
\beta^{-2n} \sum_{m=1}^n \theta^{-2mK_m} F_{2mK_m}(\theta)
\,=\, \beta^{-2n} \sum_{m=1}^n G_{2mK_m}(\theta).
\] 
Let $\P(\mbox{failure})$ denote
the probability that exact reconstruction of $f$ by $\ell^1$-minimization 
fails. By Lemma \ref{lem_key}, (\ref{Pestim30}) and by the union bound 
we finally obtain 
\begin{align}
\P(\rm{failure}) \,&\leq\, 
\P\left(\{\F_{TX} \mbox{ is injective}\} \cup \{\sup_{k \in T^c} |P_k| \geq 1\}\right)\notag\\
&\leq\, \sum_{k \in [-q,q]^d} \P(E_k) + \P(\|(N^{-1}N_0)^n\|_F \geq \kappa)
\,\leq\, D \beta^{-2n} \sum_{m=1}^n G_{2mK_m}(\theta) 
+ \kappa^{-2} M G_{2n}(\theta)\notag
\end{align} 
under the conditions
\begin{align}
a_1 \,&=\, a \,=\, \sum_{m=1}^n \beta^{n/K_m} < 1,\qquad 
a_2 + a_1 \,=\, 1 \qquad\mbox{ i.e. }\quad a_2 \,=\, 1-a,\notag\\
\frac{\kappa}{1-\kappa} \,&\leq\, \frac{a_2}{1+a_1} M^{-3/2} \,=\, 
\frac{1-a}{1+a} M^{-3/2},\notag
\end{align} 
see (\ref{cond_kappa00}). This proves Theorem \ref{thm_main0}.

\subsection{Proof of Theorem \ref{thm_main}} 

Recall that here $T$ is a random set modelled by (\ref{T_model}).
The completion of the proof of Theorem \ref{thm_main} will be slightly
more complicated as above because we still need to take the expectation
with respect to the set $T$ in Lemmas \ref{lem_H0} and
\ref{lem_Expect1}. Let us start with the expectation of $C(\A,T)$ defined
in (\ref{def_CAT}).

\begin{lemma}\label{lem_expCAT} For $\A \in P(2n,t)$ it holds
\begin{align}
\E[C(\A,T)]
&\leq\, \sum_{s=2}^n (\E|T|)^s
\sum_{R=0}^{\min\{s,t\}-1} D^{-R} 
\#\{\B \in U(2n,s), \operatorname{rank} M(\A,\B) = R\}.\notag
\end{align}
\end{lemma}
\begin{Proof} Using linearity of expectation we obtain
\begin{align}
\E[C(\A,T)]\,
&=\, \E\left[
\sum_{\substack{k_1,\hdots,k_{2n} \in [-q,q]^d, k_j \neq k_{j+1}}}
\prod_{j=1}^{2n} I_{\{k_j \in T\}}  
\prod_{A\in \A} \delta\left(\sum_{r\in A} (k_{r+1}-k_r)\right)\right]\notag\\
&=\, \sum_{\substack{k_1,\hdots,k_{2n} \in [-q,q]^d, k_j \neq k_{j+1}}}
\E\left[\prod_{j=1}^{2n} I_{\{k_j\in T\}}\right]
\prod_{A\in \A} \delta\left(\sum_{r\in A} (k_{r+1}-k_r)\right).
\notag
\end{align}
Hereby, $I_{\{k \in T\}}$ denotes an indicator variable which
is $1$ if and only if $k\in T$.
The expression $\E\left[\prod_{j=1}^{2n} I_{\{k_j \in T\}}\right]$ 
depends on how many different $k_j$'s there are. 
So once again partitions enter
the game. If $(k_1,\hdots,k_{2n}) \in ([-q,q]^d)^{2n}$ 
is a vector satisfying $k_j \neq k_{j+1}$ then we associate a partition
$\B = (B_1,\hdots,B_s)$ of $\{1,\hdots,2n\}$ such that $j$ and $j'$
are in the same set $B_i$ if and only if $k_j = k_{j'}$. 
Obviously, $j$ and $j+1$ must be contained in different blocks
for all $j$ due to the 
condition $k_j \neq k_{j+1}$ (once again we agree on the 
convention that $2n+1\equiv1$). In other words $\B$ has no adjacencies, i.e.,
$\B \in U(2n,s)$.
Now if $\B$ has $|\B| = s$ blocks then by the probability model (\ref{T_model})
for $T$ and stochastic independence
\begin{equation}\label{ETT}
\E\left[\prod_{i=1}^{2n} I_{\{k_i \in T\}}\right]
\,=\, \E\left[\prod_{j=1}^s I_{\{k_{B_j}\in T\}}\right]
\,=\, \prod_{j=1}^s \E[I_{\{k_{B_j}\in T\}}] \,=\, \tau^s,
\end{equation}
where (unambiguously) $k_{B_j} = k_i$ if $i \in B_j$.
We further introduce the notation
$\sigma_\B(r) = j$ if and only if $r \in B_j \in \B$. This leads to
\begin{align}
\E[C(\A,T)]\,
&=\, \sum_{s=2}^{2n} \tau^s \sum_{\B \in U(2n,s)}
\sum_{\substack{k_1,\hdots,k_s \in [-q,q]^d\\k_i \mbox{p.w. different}}}
 \prod_{A\in \A} \delta\left(\sum_{r\in A} 
(k_{\sigma_\B(r+1)}-k_{\sigma_\B(r)})\right).\notag
\end{align}
Clearly, the expression 
$\prod_{A \in \A} \delta\left(\sum_{r\in A} (k_{\sigma_\B(r+1)} - 
k_{\sigma_\B(r)}\right)$ is $1$ if and only if
\begin{equation}\label{solve_k}
\sum_{r \in A} (k_{\sigma_\B(r)} - k_{\sigma_\B(r+1)}) \,=\, 0
\quad\mbox{for all } A \in \A
\end{equation}
and $0$ otherwise. For $j \in \{1,\hdots,s\}$ the term
$k_j$ appears $|A_i \cap B_j|$ times as $k_{\sigma_\B(r)}$ when
$r$ runs through $A_i \in \A$. Let $M = M(\A,\B)$
denote the $t\times s$ matrix whose entries are defined by 
(\ref{defAB}). Then (\ref{solve_k}) is satisfied if and only if
$(k_1,\hdots,k_s) \in ([-q,q]^d)^s$ is contained in the kernel of $M(\A,\B)$.
Thus, if the rank of $M(\A,\B)$ equals $R$ then the number of vectors
$(k_1,\hdots,k_s) \in ([-q,q]^d)^s$ for which (\ref{solve_k}) is satisfied 
can be bounded by $D^{s-R}$ where 
$D = (2q+1)^d$. (Here we even neglected the condition 
that the $k_1,\hdots,k_s$ should be pairwise
different). So finally we obtain
\begin{align}
\E[C(\A,T)]\,
&\leq\, \sum_{s=2}^n \tau^s 
\sum_{R=0}^{\min\{s,t\}-1} D^{s-R} \#\{B \in U(n,s),
\operatorname{rank} M(\A,\B) = R\}\notag\\
&=\, \sum_{s=2}^n (\E|T|)^s \sum_{R=0}^{\min\{s,t\}-1} D^{-R} \#\{B \in U(n,s),
\operatorname{rank} M(\A,\B) = R\},\notag
\end{align}
where we substituted $\E|T| = \tau D$.
\end{Proof}
Since $\E = \E_X \E_T$ by Fubini's theorem and stochastic independence of 
$T$ and $X$ the previous result yields 
together with Lemma \ref{lem_H0}
\begin{align}
&\E[\|H_0^n\|_F^2]\notag\\
&\,\leq\, \sum_{t=1}^{\min\{n,N\}} \frac{N!}{(N-t)!} \sum_{\A \in P(2n,t)}
\sum_{s=2}^n (\E|T|)^s \sum_{R=0}^{\min\{s,t\}-1} D^{-R} \#\{B \in U(n,s),
\operatorname{rank} M(\A,\B) = R\}\notag \\
&=\, \sum_{t=1}^{\min\{n,N\}}  \frac{N!}{(N-t)!} 
\sum_{s=2}^n (\E|T|)^s \sum_{R=0}^{\min\{s,t\}-1} D^{-R} Q(2n,t,s,R)
\,=\, N^{2n} W(n,N,\E|T|,D)\notag
\end{align}
by definition (\ref{def_Q}) of the numbers $Q(2n,t,s,R)$ and by
definition (\ref{def_W}) of the function $W$. Markov's inequality yields
\[
\P(\|(N^{-1}H_0)^n\|_F \geq \kappa) \,\leq\,
N^{-2n} \kappa^{-2} \E[\|H_0^n\|_F^2] \,\leq\, \kappa^{-2} W(n,N,\E|T|,D).
\]  
We remark that by the same argument as in the proof of Theorem \ref{thm_main}
$\F_{TX}$ is injective in the event $\|(N^{-1} H_0)^n\|_F < 1$.

Let us turn now to the estimation of $\P(E_k)$. From
Lemma \ref{lem_Expect1} one realizes that we need to
estimate the expected value of $B(\A,T)$ defined in (\ref{def_BAT}). 

\begin{lemma} For $\A \in P(2Km,t)$ it holds
\begin{align}
\E[B(\A,T)]\,
&\leq\, \sum_{s=1}^{2Km} (\E|T|)^s \sum_{R=0}^{\min\{s,t\}}
D^{-R} \#\{\B \in U^*(2K,m,s), \operatorname{rank} L(\A,\B) = R\}.\notag
\end{align}
\end{lemma}
\begin{Proof} As in the proof of the previous lemma we may write
\begin{align}
\E[B(\A,T)]\,
&=\,\sum_{\substack{(k_j^{(p)})
\in ([-q,q]^d)^{2Km}\\
k_{j-1}^{(p)} \neq k_{j}^{(p)}}}
\E\left[\prod_{(p,j) \in [2K] \times [m]} I_{\{k^{(p)}_j \in T\}}\right]
\prod_{A\in \A} 
\delta\left(\sum_{(r,p) \in A} (-1)^p(k_r^{(p)} -k_{r-1}^{(p)})\right).\notag
\end{align}
Once again $\E\left[\prod_{(p,j) \in [2K] \times [m]} I_{\{k^{(p)}_j \in T\}}\right]$ 
depends on how many different $k^{(p)}$'s there are. So 
if $(k_1^{(1)},\hdots,k_m^{(2K)}) \in ([-q,q]^d)^{(2Km)}$
is a vector satisfying 
\begin{equation}\label{cond_k}
k_j^{(p)} \,\neq\, k_{j-1}^{(p)}\quad \mbox{ for all }j\in [m],~ p\in[2K]
\end{equation}
then we associate a partition $\B = (B_1,\hdots,B_s)$ of $[2K] \times [m]$ 
such that
$(p,j)$ and $(p',j')$ are contained in the same block if and only 
if $k_{j}^{(p)} = k_{j'}^{(p')}$.
Obviously, $(p,j)$ and $(p,j-1)$ cannot be contained in the same block due to
the condition (\ref{cond_k}). In other words, $\B$ belongs to $U^*(2K,m,s)$.
Now, if $\B$ has $s$ blocks, i.e., there are $s$ different values of $k^{(p)}_j$, then 
\[
\E\left[\prod_{(p,j) \in [2K] \times [m]} I_{\{k^{(p)}_j \in T\}}\right] \,=\, \tau^{s}
\]
as in (\ref{ETT}). Once more, we use the 
notation $\sigma_\B(p,j) = i$ if $(p,j) \in B_i \in \B$ and
$\sigma(p,0) = 0$. (Recall that by definition $k_0^{(p)} = k_0 = k$.)
Thus, 
\begin{align}
\E[B(\A,T)] \,=\,
\sum_{s=1}^{2n}\tau^s \sum_{\B \in U^*(2K,m,s)}
\sum_{\substack{ k_1,\hdots,k_s \in [-q,q]^d\\ k_i \displaystyle{\mbox{ p.w. different}}}}
\prod_{A \in \A}\delta\left(\sum_{(p,j) \in A} (-1)^p (k_{\sigma_\B(p,j)} - k_{\sigma_\B(p,j-1)})\right).\notag
\end{align}
The term
$\prod_{A \in \A}\delta\left(\sum_{(p,j) \in A} (-1)^p (k_{\sigma_\B(p,j)} - k_{\sigma_\B(p,j-1)})\right)$
contributes to the sum if and only if
\[
\sum_{(p,j) \in A} (-1)^p (k_{\sigma_\B(p,j)} - k_{\sigma_\B(p,j-1)}) 
\,=\, 0 \quad \mbox{ for all } A\in \A.
\]
By definition (\ref{def_LAB}) of the matrix $L(\A,\B)$ and since $k_0 = k$ 
this is equivalent to
\begin{equation}\label{LAB_eq}
L(\A,\B) (k_1,\hdots,k_s)^T \,=\, k v(\A,\B),
\end{equation}
where $v=v(\A,\B)$ is the $t$-dimensional vector with entries
\[
v_i \, = \, \sum_{(p,1)\in A_i} (-1)^p, \quad i=1,\hdots, t.
\]
(If $d>1$ then (\ref{LAB_eq}) has to interpreted vector-valued, i.e., for each
component of $k \in [-q,q]^d$ and of $k_1,\hdots,k_s \in [-q,q]^d$ we have one equation
with the same $L(\A,\B)$ and the same $v(\A,\B)$.)
If the rank of $L(\A,\B)$ equals $R$ then we can bound the number 
of solutions to (\ref{LAB_eq})
by $D^{s-R}$. Hence, we obtain the bound
\[
\E[B(\A,T)] \,\leq\, \sum_{s=1}^{2Km} \tau^s \sum_{R=0}^{\min\{s,t\}} D^{s-R} 
\#\{\B \in U^*(2K,m,s), \operatorname{rank} L(\A,\B) = R\}.
\]  
Since $\E|T| = \tau D$ this proves the lemma.
\end{Proof}
Together with Lemma \ref{lem_Expect1} the previous result yields
\begin{align}
\E\left[|((HR_T)^m \sigma)_{k}|^{2K}\right]
\,&\leq\,  \sum_{t=1}^{\min\{Km,N\}} \frac{N!}{(N-t)!} 
\sum_{s=1}^{2Km} (\E|T|)^s \sum_{R=0}^{\min\{s,t\}} Q^*(2K,m,t,s,R)\,D^{-R}\notag\\
&=\, N^{2Km} Z(K,m,N,\E|T|,D)\notag
\end{align}
where $Q^*(2K,m,t,s,R)$ are the numbers defined in (\ref{def_Qs}). 
By (\ref{prob_Ek}) we obtain
\[
\P(E_k) \,\leq\, \beta^{-2n} \sum_{m=1}^n Z(K_m,m,N,\E|T|,D)
\]
Finally, let $\P(\rm{failure})$ denote
the probability that exact reconstruction of $f$ fails.
By Lemma \ref{lem_key}, (\ref{Pestim_3}), (\ref{prob_ET}) 
and using that 
$\{\F_{TX} \mbox{ is not injective}\} \subset \{\|(N^{-1}H_0)^n\|_F \geq \kappa\}$
we obtain
\begin{align}
\P(\rm{failure}) \,&\leq\, 
\P\left(\{\F_{TX} \mbox{ is injective}\} \cup \{\sup_{k \in T^c} |P_k| \geq 1\}\right)\notag\\
&\leq\, \sum_{k \in [-q,q]^d} \P(E_k) + \P(\|(N^{-1}N_0)^n\|_F \geq \kappa) + \P(|T| \geq (\alpha + 1)\E|T|)\notag\\
&\leq\, D \beta^{-2n} \sum_{m=1}^n Z(K_m,m,N,\E|T|,D) + \kappa^{-2} W(n,N,\E|T|,D)
+ \exp(- \frac{3\alpha^2}{6+2\alpha} \E|T|)\notag
\end{align} 
under the conditions
\begin{align}
a_1 \,&=\, a \,=\, \sum_{m=1}^n \beta^{n/K_m} < 1,\qquad
a_2 + a_1 \,=\, 1 \qquad\mbox{ i.e. }\quad a_2 \,=\, 1-a,\notag\\
\frac{\kappa}{1-\kappa} \,&\leq\, \frac{a_2}{1+a_1} ((\alpha+1)\E|T|)^{-3/2} \,=\, 
\frac{1-a}{1+a}((\alpha+1)\E|T|)^{-3/2},\notag
\end{align} 
see (\ref{cond_kappa}). This proves Theorem \ref{thm_main}.

\subsection{Proof of Corollary \ref{cor1}}

We have to show that a finer analysis of the 
probability bound (\ref{Pbound0}) of Theorem \ref{thm_main0}
gives Corollary \ref{cor1}.
We first claim that the associated Stirling numbers satisfy the estimate
\begin{equation}\label{S2_claim}
S_2(n,k) \,\leq \, (3n/2)^{n-k} \quad \mbox{for all } k=1,\hdots,\lfloor n/2 \rfloor.
\end{equation}
Indeed, the claim is true for $S_2(1,k)=0$ and $S_2(2,1) = 1$. Now suppose,
the claim is true for all $S_2(m,k)$ with $m < n$. Then from the recursion formula
(\ref{rec_formula}) it follows
\begin{align}
S_2(n,k) \,&=\, k S_2(n-1,k) + (n-1) S_2(n-2,k-1) \notag\\
\,&\leq\, k (3(n-1)/2)^{n-k-1} + (n-1)(3n/2-3)^{n-k-1} 
\,\leq\, (n-1+k)(3n/2)^{n-k-1}\notag\\ 
&\leq\, (3n/2)^{n-k}\notag
\end{align}
since $n-1+k \leq 3n/2$. This proves (\ref{S2_claim}). Pluggin this into the definition
of $G_{2n}$ yields
\begin{align}
G_{2n}(\theta) \,&=\, \theta^{-2n} \sum_{k=1}^n S_2(2n,k) \theta^{k} 
\,\leq\, \theta^{-2n} \sum_{k=1}^n (3n)^{2n-k} \theta^k
\,=\, (3n/\theta)^{2n} \sum_{k=1}^n (\theta/3n)^k\notag\\
&=\, (3n/\theta)^{2n} \frac{(\theta/3n)^{n+1}-(\theta/3n)}{(\theta/3n)-1}
\,=\, (3n/\theta)^{2n-1}  \frac{(\theta/3n)^{n}-1}{(\theta/3n)-1}.\notag
\end{align}
Now assume we have chosen $n$ such that $n\leq \theta/6$. Then we further obtain
\[
G_{2n}(\theta) \,\leq\, (3n/\theta)^{n-1}.
\] 
Now consider the term $D\beta^{-2n} \sum_{m=1}^n G_{2mK_m}(\theta)$ from the 
probability bound (\ref{Pbound0}). We choose $K_m=r(n/m)$ where $r$ denotes
the function that rounds to the nearest integer. Then it is easy to see that 
\[
m K_m \in \{\lceil 2n/3 \rceil,\hdots,\lfloor 4n/3\rfloor\},
\quad m \in \{1,\hdots,n\}.
\]
Thus,
\[
\sum_{m=1}^n G_2mK_m(\theta) \,\leq\, 
n \max_{k \in \{\lceil 2n/3 \rceil,\hdots, \lfloor 4n/3 \rfloor\}} G_{2k}(\theta)
\,\leq\, n \left(\frac{4n}{\theta}\right)^{2n/3 -1}
\]
provided $k \leq \theta/6$ for all 
$k \in \{\lceil 2n/3 \rceil,\hdots,\lfloor 4n/3\rfloor\}$, i.e.,
$6 \lfloor 4n/3 \rfloor \leq \theta$. 
This yields
\[
D\beta^{-2n} \sum_{m=1}^n G_{2mK_m}(\theta) \,\leq\, 
D n \left(\frac{4n}{\theta}\right)^{-1} 
\left(\beta^{-3} \frac{4n}{\theta}\right)^{2n/3}.
\]
In order to make this expression small it is certainly a good strategy to make
the last term smaller than $1$. Indeed, choose
\begin{equation}\label{choice_n}
n \,=\, n(\theta) \,:=\, \left\lfloor \frac{\beta^3 \theta}{8} \right\rfloor
\end{equation}
implying $\beta^{-3} 4n/\theta \leq 1/2$. (This choice for $n$
is certainly valid since $\beta < 1$ as it must satisfy 
condition (\ref{cond_kappa_0}).) We obtain
\[
D \beta^{-2n} \sum_{m=1}^n G_{2mK_m}(\theta)
\,\leq\, \frac{1}{4} D \theta\, 2^{-2n(\theta)/3}.
\]
A simple calculation yields that the latter term is less than $\epsilon/2$ if
\[
\frac{2\ln(2)}{3} n(\theta) - \ln(\theta) \geq \ln(D) + \ln(\epsilon^{-1}) - \ln(2).
\] 
Furthermore, a simple numerical test shows that a valid choice for $\beta$
is $\beta = 0.47$. The corresponding $a = \sum_{m=1}^n \beta^{n/K_m}$ is always 
less than $0.957$ and $n(\theta) \approx \lfloor 0.013\, \theta \rfloor$.
Recalling that $\theta = M/N$ it follows that there exists a constant $C_1$
such that $D \beta^{-2n} \sum_{m=1}^n G_{2mK_m} (\theta) \leq \epsilon/2$
provided
\[
N \,\geq\, C_1 M (\ln(D) + \ln(\epsilon^{-1})).
\]
Now consider the other term $M\kappa^{-2} G_{2n}(\theta)$ in the probability bound 
(\ref{Pbound0}). We choose $\kappa$ such that there is equality in 
(\ref{cond_kappa_0}), i.e.,
\[
\kappa \,=\, \frac{(1-a)/(1-a)M^{-3/2}}{1+(1-a)/(1+a)M^{-3/2}} 
\,\geq\, \frac{1-a}{2(1+a)} M^{-3/2}.
\]  
Hence,
\[
M\kappa^{-2} G_{2n}(\theta) \,\leq\, \left(\frac{1-a}{2(1+a)}\right)^2 M^4 
G_{2n}(\theta) 
\]
Now we do not have the freedom anymore to choose $n$. We have to make the same 
choice (\ref{choice_n}) as above. This yields
\[
M\kappa^{-2} G_{2n(\theta)}(\theta)
\,\leq\,  \left(\frac{2(1+a)}{(1-a)}\right)^2 M^4 
\left(\frac{3\beta^3}{8}\right)^{n(\theta)-1}.
\]
Requiring that the latter expression is less than $\epsilon/2$ is equivalent
to
\[
(n(\theta)-1)\ln\left(\frac{8}{3\beta^3}\right) \,\geq\, 
\ln\left(8\left(\frac{1+a}{1-a}\right)^2\right) + 4 \ln(M) + \ln(\epsilon^{-1}).
\]
As already remarked the choice $\beta=0.47$ results in $a \leq 0.957$ and
$n(\theta) \approx \lfloor 0.013\,\theta\rfloor$.
Hence, $\ln(8/(3\beta^3)) \approx 3.2459$ and $\ln(8((1+a)/(1-a))^2) \approx 9.7153$.
Since $M\leq D$ there exists a constant $C_2$ 
(whose precise value may be calculated
from the numbers above) such that
$M\kappa^{-2} G_{2n(\theta)}(\theta) \leq \epsilon/2$
provided
\[
N \geq C_2 M (\ln(D) + \ln(\epsilon^{-1})).
\]
Choosing $C:= \max\{C_1,C_2\}$ completes the proof of Corollary \ref{cor1}.

We remark that analyzing numerical plots 
for $\beta^{-2n'(\theta)} \sum_{m=1}^{n'(\theta)} G_{2mK_m(\theta)}(\theta)$ and 
$G_{2n'(\theta)}(\theta)$ for $n'(\theta) = \lceil \theta/12 \rceil$ indicates
that one may choose the constant $C$ much smaller as the ones resulting
from the theoretical analysis above. 
It seems that $C \lessapprox 20$ is a valid choice.

\subsection{Remarks}
\label{sec_remarks}

We conclude this section with some remarks.
\begin{itemize}
\item[(a)] Let us give a more detailed reason why we believe that the 
probilistic model for the ``sparsity set'' $T$ is likely to 
give better probability bounds for exact reconstruction than the deterministic
approach holding for all $T$ of a given size. Indeed the main difference
in the two previous proofs lies in the estimation of $C(\A,T)$ and $B(\A,T)$ 
defined in (\ref{def_CAT}) and (\ref{def_BAT}).
If $|\A| = t$ then
for deterministic $T$ we used the estimation (\ref{CAT_deterministic}), i.e., 
$C(\A,T) \,\leq\, |T|^{2n-t+1}$. Indeed, if $T$ is an arithmetic progression then
$C(\A,T)$ may come very close to this upper bound. However, for generic sets $T$
the bound is quite pessimistic.
In fact, in the probabilistic model the expected size of $C(\A,T)$ can be bounded
by
\begin{equation}\label{Exp_CAT}
\E[C(\A,T)] \,\leq\, \sum_{s=2}^n (\E|T|)^s
\sum_{R=0}^{\min\{s,t\}-1} D^{-R} 
\#\{\B \in U(2n,s), \operatorname{rank} M(\A,\B) = R\},\notag
\end{equation}
see Lemma \ref{lem_expCAT}.
In particular, if $D$ is large (and $\E|T|$ not too small) then 
the latter estimate should be much better. Let us illustrate this with two
examples. 
\begin{enumerate}
\item Let $\A=\{\{1,2,3,5\},\{4,6\}\}$, i.e., $2n = 6$ and $t=2$. Then
(\ref{CAT_deterministic}) yields $C(\A,T) \leq |T|^5$ while computing 
\ref{Exp_CAT}) explicitly gives
\[
\E[C(\A,T)] \,=\, D^{-1}\left[(\E|T|)^2 + 10 (\E|T|)^3 + 20 (\E|T|)^4
+9 (\E|T|)^5 + (\E|T|)^6\right].
\]
Clearly, if $D$ is sufficiently large then the probabilistic estimate
is much better than the deterministic one.
\item Let $\A=\{\{1,2,3\},\{4,5,6\}\}$, so again $2n=6$ and $t=2$. Then
the deterministic estimate gives again $C(\A,T) \leq |T|^5$ 
while (\ref{Exp_CAT}) results in
\begin{align}
\E[C(\A,T)] \,&\leq\, (\E|T|)^2 + 3 (\E|T|)^3 + (\E|T|)^4 \notag\\
& + D^{-1}\left[7 (\E|T|)^3 +19 (\E|T|)^4 + 9 (\E|T|)^5 + (\E|T|)^6\right].\notag
\end{align}
So here one has to choose both $\E|T|$ and $D >> \E|T|$ large to see
that potentially the probabilistic estimate is much better.
\end{enumerate}
\item[(b)] {\bf Discrete Fourier transforms:} 
The whole proofs work without essential change if one replaces
our setting by the following one similar to the situation investigated
by Candes, Romberg and Tao in \cite{CRT1}. 
Consider functions on the cyclic group $\Z_p^d = \{0,\hdots,p-1\}^d$, $p \in \N$, 
rather than on $[0,2\pi]^d$.
The discrete Fourier transform is defined by 
\[
\hat{f}(\omega) \,:=\, \sum_{x \in \Z_p^d} f(x) e^{-2\pi i x\cdot \omega/p},\quad
\omega \in \Z_p^d.
\]
We draw $x_1,\hdots,x_N$ from the uniform distribution
on $\Z_p^d$. Note that in contrast to sampling from $[0,2\pi]^d$
it may occur with non-zero probability that some
elements of $\Z_p^d$ are drawn more than
once. But this will not do much harm. 
 
Let $f$ be such that $\hat{f}$ is a sparse vector on $\Z_p^d$. Once again
we try to reconstruct $f$ from its sample values $f(x_j)$ by minimizing
the $\ell^1$-norm of $\hat{f}$ under the constraint that the observed values $f(x_j)$ 
are matched. 

Theorems \ref{thm_main0} and \ref{thm_main} will also apply to this situation.
Indeed, the only thing that differs in the proofs 
is that we have to calculate modulo 
$p$ in the definition of $C(\A,T)$ and $B(\A,T)$, see (\ref{def_CAT}) 
and (\ref{def_BAT}). This is apparent from (\ref{Expect_Trig}) where
the integral is replaced by a sum of exponentials.
Nevertheless, the deterministic and probabilistic 
estimates for the quantities $C(\A,T)$ and $B(\A,T)$ still hold and so
everything goes through in completely the same manner.

Of course, one can also exchange the role of $f$ and $\hat{f}$,
aiming at reconstructing a sparse signal on $\Z_p^d$ from
random samples of its Fourier transform. Indeed, this situation
is investigated in \cite{CRT1} with a different probability model
for the sampling points. In other words, we presented a slightly different 
approach for the main result in \cite{CRT1}.
\end{itemize}

\section{Some more on set partitions}
\label{sec_Qn}

From Theorem \ref{thm_main} we realize that we have to investigate
the functions $F_n(\theta)$ connected to set partitions in $P(n,t)$ and
also the numbers $Q(n,t,s,R)$ and $Q^*(K,m,t,s,R)$, respectively. 
We already gave some information
on the number $S_2(n,t)$ of partitions in $P(n,t)$ earlier.
Let us be a bit more detailed here. Clearly, by definition (\ref{def_Fn})
of $F_n$ and the generating function (\ref{expgen_func})
we see that
\[
F_{n}(\theta) \,=\, \sum_{k=1}^{\lfloor n/2 \rfloor} S_2(n,k) \theta^k.
\]   
(This follows also directly from the proof of Theorem \ref{thm_main0}.)
In particular, $F_{2n}$ is a polynomial of degree $n$. There are different
ways of computing $F_{n}$ explicitly. One possibility is to use the generating
function (\ref{expgen_func}) leading to
\[
F_n(\theta) \,=\, 
\frac{\partial^n}{\partial x^n} \exp(\theta(e^x-x-1))_{\vert_{x=0}}\,.
\]
One may also compute the numbers $S_2(n,k)$ explicitly. Indeed, 
differentiating (\ref{expgen_func}) $k$ times with respect to $y$ and
setting $y=0$ yields
\[
\sum_{n=1}^\infty S_2(n,k) \frac{x^n}{n!} \,=\, 
\frac{1}{k!}\left(e^x-x-1\right)^k.
\]
Expanding the right hand side into a power series and comparing coefficients 
yields (after some computations)
\begin{equation}\label{S2explicit}
S_2(n,k) \,=\, \frac{1}{k!} \left(k^n + 
\sum_{j=1}^{k-1} (-1)^j\left(\begin{matrix}k\\ j\end{matrix}\right) 
\sum_{\ell=0}^j \left(\begin{matrix} j\\ \ell\end{matrix}\right)
\frac{n!}{(n-\ell)!}(k-j)^{n-\ell}\right)
\end{equation}
valid for $n \geq 2k$ (otherwise $S_2(n,k) = 0$).
In the special case $k=2$ we obtain
$S_2(n,2) = 2^{n-1}-n-1$.
Further, a combinatorial argument shows that $S_2(2n,n) = \frac{2n!}{2^n n!}$. (One uses that $P(2n,n)$ consists
only of partitions where each block has precisely $2$ elements.)

Let us give the first of the functions $F_{2n}$ explicitly in the
following list,
\begin{align}
F_{2}(y) \,&=\, y,\qquad 
F_{4}(y) \,=\, y + 3 y^2, \qquad
F_{6}(y) \,=\, y + 25 y^2 + 15 y^3,\notag\\
F_{8}(y) \,&=\, y + 119 y^2 + 490 y^3 + 105 y^4,\qquad
F_{10}(y)\,=\, y + 501 y^2 + 6825 y^3 + 9450 y^4 + 
945 y^5,\notag\\
F_{12}(y) \,&=\, y + 2035 y^2 + 74316 y^3 + 302995y^4
+190575y^5 + 10395y^6\notag. 
\end{align}
Of course, explicit values of $S_2(n,k)$ can be read off this list.

Now consider the number $p_n = \sum_k S_2(n,k)$ of
all partitions of $[n]$ into subsets having at least two elements.
Setting $y=1$ in the exponential generating 
function (\ref{expgen_func}) yields
\begin{equation}\label{gen_func_pn}
\sum_{n=1}^\infty p_n \frac{x^n}{n!} \,=\, \exp(e^x-x-1).
\end{equation}

Unfortunately, much less is known
about the number of partitions in $U(n,s)$. As already mentioned, it was only very recently
that D. Knuth \cite{Knuth} posed the problem of determining $|U(n,s)|$.
Let us denote by $u_n = \sum_{k=2}^n |U(n,k)|$
the number of all partitions of $\{1,\hdots,n\}$ having no adjacencies
(recall that $U(n,1) = \emptyset$).
Recently, it was proved in \cite{Callan}
that $u_n = p_n$. So (\ref{gen_func_pn}) is also the exponential generating 
function of the numbers $u_n$.
Concerning the size of $U^*(K,m,s)$, up to now, we cannot say more than 
that it is bounded
by the number of all partitions into $s$ blocks of a set with $Km$ elements, i.e.,
by the (ordinary) Stirling number of the second kind $S(Km,s)$. If $m=1$ then $|U^*(K,1,s)| = S(K,s)$
as already remarked.
The Stirling numbers $S(n,k)$ have the generating function 
\cite{Riordan,Stanley}
\begin{equation}\label{expgen_func_Stirling}
\sum_{n=1}^\infty \sum_{k=1}^n S(n,k) y^k \frac{x^n}{n!}
\,=\, \exp(y(e^x-1)).
\end{equation}
Let us denote $u^*_{K,m} = \sum_{k=1}^{Km} |U^*(K,n,k)|$. Then clearly 
$u^*_{K,m} \leq b_{Km} = \sum_{k=1}^n S(Km,k)$ with equality if $m=1$. 
A lower bound for $u^*_{K,m}$ is given by the numbers $p_{Km}$.

Now some elementary observations concerning the numbers $Q(n,t,s,R)$ and 
$Q^*(K,n,t,s,R)$ can be made.
Disregarding the rank of $M(\A,\B)$, the number of all pairs $(\A,\B)$ with
$\A \in P(n,t)$ and $\B \in U(n,s)$ is $|P(n,t)| \times |U(n,s)|$, hence, 
$\sum_{R=0}^{\min\{s,t\}} Q(n,t,s,R) = |P(n,t)| \times |U(n,s)|$
and similarly for $Q^*(K,m,t,s,R)$. Summing also over $t$ and $s$ gives
\[
\sum_t \sum_s \sum_R Q(n,t,s,R) \,=\, u_n p_n \,=\, p_n^2
\]
and $\sum_{t,s,R} Q^*(K,m,t,s,R) = p_{Km} u^*_{K,m}$. In the following
table we give some values of $p_n$, and $b_n$ for even $n=2,4,6,\hdots$ 
(we omit the odd
numbers since we do not need them for Theorem \ref{thm_main}).

\bigskip

\noindent
\begin{tabular}{|l|c|c|c|c|c|c|c|c|} \hline
$n$  &  2  &  4  &  6  &  8   &  10  &  12  &  14 & 16 \\ \hline 
$p_n = u_n$ &  1  &  4  &  41 &  715 &  17 722 & 580 317 & 24 011 157 & 1 216 070 380\\ \hline 
$b_n = u^*_{n,1}$ & 2 & 15 & 203 & 4140 & 115 975 & 4 213 597 &  190 899 322 
& 10 480 142 147\\ \hline
\end{tabular}

\bigskip

We determined $Q(n,t,s,R)$ and $Q^*(K,m,t,s,R)$ for certain small $n,K,m$ 
on a computer in the following way.
First all partitions in $P(n,t)$ and $U(n,s)$ (resp. $U^*(K,m,s)$) 
are computed recursively. For $P(n,t)$ we have the following procedure:
\begin{enumerate}\itemsep-0.5pt
\item if $n < 2$ or $t > n/2$ then RETURN $P(n,t) = \emptyset$. 
\item if $t = 1$ then RETURN $P(n,1) = \{\{1,\hdots,n\}\}$.
\item $P(n,t) = \emptyset$
\item compute (recursively) $P(n-1,t)$ and $P(n-2,t-1)$.
\item for each $\A \in P(n-1,t)$: \\
\phantom{for } for $j$ from $1$ to $t$: \\
\phantom{for for } create new partition $\A'$ by adding 
          the element $n$ to the $j$-th subset of $\A$\\
\phantom{for for } add $\A'$ to $P(n,t)$ 
\item for each $\A \in P(n-2,t-1)$:\\
\phantom{for } for $\ell$ from $1$ to $n-1$:\\
\phantom{for for } create new partition $\A'$ from $\A$ 
by incrementing each
element $i \in [n-2]$ by $1$ if $i \geq \ell$\\ 
\phantom{for for bb} and adding the subset $\{\ell,n\}$\\
\phantom{for for } add $\A'$ to $P(n,t)$ 
\item RETURN $P(n,t)$
\end{enumerate}
We remark that from this procedure also 
the recursion formula (\ref{rec_formula}) follows.

The partitions in $U(n,s)$ are determined by first computing the set $V(n,s)$ 
of {\it all} partitions of $[n]$ into $s$ blocks and then omitting those that have
adjacencies. Similarly $U^*(K,n,s)$ is computed. Hereby, we have the following
recursive procedure to compute $V(n,s)$:
\begin{enumerate}\itemsep-0.5pt
\item if $s=1$ RETURN $V(n,s) = \{\{1,\hdots,n\}\}$
\item if $s = n$ RETURN $V(n,s) = \{\{1\},\hdots,\{n\}\}$
\item $V(n,s) = \emptyset$ 
\item compute (recursively) $V(n-1,s)$ and $V(n-1,s-1)$
\item for each $\A \in V(n-1,s)$:\\
\phantom{for } for $j$ from $1$ to $s$:\\
\phantom{for for } create new partition $\A'$ by adding 
          the element $n$ to the $j$-th subset of $\A$\\
\phantom{for for } add $\A'$ to $V(n,s)$ 
\item for each $\A \in V(n-1,s-1)$:\\
\phantom{for } create new partition $\A'= \A \cup \{n\}$\\ 
\phantom{for } add $\A'$ to $V(n,s)$
\item RETURN $V(n,s)$
\end{enumerate}
One may easily deduce the recursion formula $S(n,k) = kS(n-1,k) + S(n-1,k-1)$
for the Stirling numbers of the second kind $S(n,k) = |V(n,k)|$ from this
procedure. 

After determining $P(n,t)$ and $U(n,s)$ 
for each pair $(\A,\B)$ with $\A \in P(n,t)$ and $\B \in U(n,s)$ 
(or $\B \in U^*(K,n,s)$ resp.) we set up the matrix $M(\A,\B)$ (or $L(\A,\B)$),
see (\ref{defAB}) and (\ref{def_LAB}), and compute its rank. By counting
the number of matrices $M(\A,\B)$ that have rank $R$ we 
determine $Q(n,t,s,R)$ or
$Q^*(K,m,t,s,R)$, respectively. The results of these computations for 
certain $n,K,m$ are given in the appendix.
Considering the table of the numbers $p_n$ (recall that
$p_n^2$ equals the overall number of matrices whose rank has to be determined)
we realize that this procedure is practicable only for small values of $n$.
Even for $n=10$ the computing time reaches several days and for $n=14$ it seems
impossible to do the task in a reasonable time as 
$p_{14}^2 = 576 535 660 478 649$.

The following lemma is concerned with $Q(n,t,s,0)$ for some special
cases.

\begin{lemma}\label{lem_Qn0}
\begin{itemize}
\item[(a)] $Q(n,1,s,0) \,=\, |U(n,s)|$.
\item[(b)] It holds
$
Q(2n,2,2,0) 
= \frac{(2n)!}{2(n!)^2} - 1
$ 
and $Q(2n,2,2,1) = 2^{2n-1}-2n-\frac{(2n)!}{2(n!)^2}$.
\item[(c)] If $n > 2$ and $2s \geq 3n$ then
$Q(2n,n,s,0) = 0$.
\item[(d)] If $t \neq 1$ then
$Q(2n,t,2n,0) = 0$.
\item[(e)] If $n> 3$ and $3t\geq 2n$ then $Q(2n,t,2n-1,0) = 0$.
\end{itemize}
\end{lemma}
\begin{Proof} 
(a) There is only one partition in $P(n,1)$ and the maximal rank of $M(\A,\B)$ is 
$\min\{t,s\}-1 = 0$. Thus, $Q(n,1,s,0) = |U(n,s)|$.

(b) Clearly, $U(2n,2)$ consists of only $1$ partition $\B = (B_1,B_2)$, i.e.,
\begin{align}
B_1 \,=\, \{1,3,5,\hdots,2n-1\},\qquad
B_2 \,=\, \{2,4,6,\hdots,2n\}.\notag
\end{align}
The associated matrix $M=M(\A,\B)$, $\A=\{A_1,A_2\} \in P(2n,2)$ has entries 
\[
M_{i,j} \,=\, |A_i \cap B_j| - |(A_i + 1) \cap B_j| \,=\, |A_i \cap B_j| - |A_i \cap (B_j -1)|
\,=\, |A_i \cap B_j| - |A_i \cap B_{3-j}|
\]
since $B_1-1=B_2$ and $B_2-1=B_1$. Thus, $M(\A,\B)$ has rank $0$, i.e.,
$M(\A,\B) = 0$ if and only if
\begin{equation}\label{cond22}
|A_1 \cap B_1| \,=\, |A_1 \cap B_2| 
\quad \mbox{and} \quad |A_2 \cap B_1| \,=\, |A_2 \cap B_2|.
\end{equation}
So $A_1$ and $A_2$ 
must have the same number of elements from $B_1$ and from $B_2$.
So we can construct all possible partitions $\A$ satisfying (\ref{cond22})
in the following way. Choose $m \in \{1,\hdots,n-1\}$ and then form $A_1$ by taking
$m$ elements from $B_1$ and $m$ elements from $B_2$. The set $A_2$ is formed of
all the remaining elements. Then (\ref{cond22}) is clearly satisfied.
We can do this in $\left(\begin{matrix} n\\m\end{matrix}\right)^2$ different ways. 
However, if we run with $m$ 
through $\{1,\hdots,n-1\}$ every possible partition appears once as $\{A_1,A_2\}$ and
once as $\{A_2,A_1\}$, so that altogether we have the formula
\[
Q(2n,2,2,0) \,=\, 2^{-1} \sum_{m=1}^{n-1} \left(\begin{matrix} n\\ m \end{matrix}\right)^2
\,=\, \frac{(2n)!}{2(n!)^2} - 1.
\]
The second equality follows from the fact that 
$\sum_{m=0}^n  \left(\begin{matrix} n\\ m \end{matrix}\right)^2 =  
\left(\begin{matrix} 2n\\ n \end{matrix}\right)$, see e.g. \cite{Weisstein}.
Now the second assertion follows easily since 
\[
Q(2n,2,2,1) \,=\, |P(2n,2)| - Q(2n,2,2,0)
\,=\, 2^{2n-1}-1-2n- Q(2n,2,2,0).
\]

(c)-(e) For all the remaining cases we have to prove that for all relevant
partitions $\A \in P(2n,t), \B \in U(2n,s)$ we never have $M(\A,\B) = 0$
(the zero-matrix). Observe that $M(\A,\B) = 0$ means that
\begin{equation}\label{cond_sets}
|A \cap B| \,=\, |(A+1) \cap B| \quad \mbox{for all } A \in \A, B \in \B 
\end{equation}
(where $A+1$ is computed modulo $n$ as usual). 
So for all three cases we assume that $\A \in P(2n,t)$ and $\B \in U(2n,s)$
are given (with $t,s$ satisfying the respective conditions)
and show that the condition (\ref{cond_sets}) leads to a contradiction. 

(c) Clearly, a partition $\A$ in $P(2n,n)$ has 
only subsets consisting of precisely $2$ elements. The condition
$2s \geq 3n$ implies that a partition $\B \in U(2n,s)$ has
at least $n$ singletons (i.e. subsets consisting of only one element).
Indeed, if there would be less than $n$ singletons than the overall
number of elements would be larger than $n-1 + 2(s-(n-1))$ (i.e. $n-1$ sets
with $1$ element and $(s-(n-1))$ sets with at least $2$ elements). 
Since $n-1 + 2(s-(n-1)) = 2s-n+1 \geq 3n-n+1 = 2n+1$ 
this produces a contradiction
as there are only $2n$ elements. 

Now, if $\{k\}$ is a singleton of $\B$ and $k \in A$ for $A \in \A$ 
then condition (\ref{cond_sets}) implies that also $(k-1) \in A$. 
As all subsets $A$ in $\A$ have precisely two elements this means
that $A = \{k-1,k\}$. Using once more (\ref{cond_sets}) we further
see that this implies that neither $\{k-1\}$ nor $\{k+1\}$ can be singletons
in $\B$. So $\A$ has the form
\[
\{ \{1,2\},\{3,4\},\hdots,\{2n-1,2n\}\}
\]
and the singletons of $\B$ are $\{2\},\{4\},\hdots,\{2n\}$ up to shifting
all elements by $1$ (modulo $n$). We still have to distribute the remaining
numbers $1,3,5,\hdots,2n-1$ onto subsets in $\B$. If $1 \in B\in \B$ then condition
(\ref{cond_sets}) with $A=\{1,2\}$ tells us that also $3 \in B$. The same
argument for $3$ and $A = \{3,4\}$ implies that also $5 \in B$ and so on.
So $B=\{1,3,5,\hdots,2n-1\}$ and thus, $s=n+1$. Since $n > 2$ 
this is a contradiction
to $s \geq 3n/2$. Thus, there is no pair of partitions 
$\A \in P(2n,n), \B \in U(2n,s)$ with $M(\A,\B) = 0$.

(d) The only partition in $U(2n,2n)$ is $\B = \{\{1\},\{2\},\hdots, \{2n\}\}$.
Thus the condition (\ref{cond_sets}) implies that whenever $j \in A \in \A$
then also $j-1 \in A$. As $j$ is arbitrary this means that the only
possibility for $A$ is $\{1,2,\hdots,2n\}$, i.e., $t=1$. 

(e) The condition on $t$ implies that there is at least one subset 
$A_1 \in \A \in P(2n,t)$ that has precisely $2$ elements. Moreover,
any partition in $U(2n,2n-1)$ has precisely $2n-2$ singletons and
one subset $B_1$ consisting of precisely $2$ elements. We write
$A_1 = \{j_1,j_2\}$, $j_2 > j_1$ and $B_1=\{k_1,k_2\}$. 

We distinguish two cases.
Let us first assume $j_2 \neq j_1+1$ and $\{j_1,j_2\} \neq \{1,2n\}$.
All singletons in $\B$ are given
by $\{k\}$ with $k\neq k_1,k_2$. Checking the condition (\ref{cond_sets})
with $A_1$ and $\{k\}$ shows that necessarily $j_2 \neq k$ for all $k\neq k_1,k_2$.
Without loss of generality this means $j_2 = k_2$. Condition (\ref{cond_sets})
with $A_1$ and $B_1$ thus yields
\[
|\{j_1,j_2\} \cap \{k_1,j_2\}| = |\{j_1+1,j_2+1\} \cap \{k_1,j_2\}|
\]
It is not possible that the sets on both sides have both cardinality $2$. Thus,
the relation implies $j_1 \neq k_1$. Moreover, since by assumption $j_1 +1 \neq j_2$
either $k_1 = j_1+1$ or $k_1 = j_2 + 1$. 
In both cases the singleton $\{j_1\}$ belongs
to $\B$. Condition (\ref{cond_sets}) yields 
$|\{j_1,j_2\} \cap \{j_1\}| = |\{j_1+1,j_2+1\} \cap \{j_1\}|$ which
is not possible since $j_1 \neq j_2+1$ by the assumptions $j_2 \geq j_1$
and $\{j_1,j_2\} \neq \{1,2n\}$.

Next we treat the case $A_1 = \{j_1,j_2\} = \{j,j+1\}$. Without loss
of generality we may assume $j=1$, so $A_1 = \{1,2\}$.
Checking condition (\ref{cond_sets}) with $A_1$ and $\{k\}$, $k \neq 1,2$
shows that $k \neq 1,3$. Thus $B_1 = \{1,3\}$ and the singletons
of $\B$ are the sets $\{2\},\{4\},\{5\},\{6\},\hdots,\{2n\}$.
Then condition (\ref{cond_sets}) with $A_1$ and $B_1$ is satisfied.
Now, let $A$ be the subset of $\A$ containing the element $3$ and 
write $A = \{3\} \cup A'$. Then condition (\ref{cond_sets}) with $B= \{4\}$
reads $|(A' \cup \{3\}) \cap \{4\}| = |((A'+1) \cup \{4\}) \cap \{4\}| = 1$.
Thus, $A'$ and hence also $A$ must contain the element $4$. We may continue 
in this way to show that $A=\{3,4,5,\hdots,2n\}$. In particular, $t=2$.
Since $n > 3$ this is a contradiction to $3t > 2n$.
\end{Proof}
One may compare the assertions of this lemma with the tables in the
appendix.
For $Q^*(K,m,t,s,0)$ certainly a similar analysis can be done but
we have not further pursued this issue here.

\section{Bounds for the probability of exact reconstruction}
\label{sec_bound}

In this section
we illustrate the bounds in Theorems \ref{thm_main0} and \ref{thm_main} 
for the probability of exact reconstruction by drawing 
some plots. Hereby, we always plotted the bound of the probability
of {\it failure} of exact reconstruction, i.e., $1$ minus the
expressions in (\ref{Pbound0}) and (\ref{Pbound}).

\begin{figure}\label{fig_prob_thm0}
\begin{center}
\parbox[t]{12cm}{
\includegraphics[width=12cm,height=10cm]{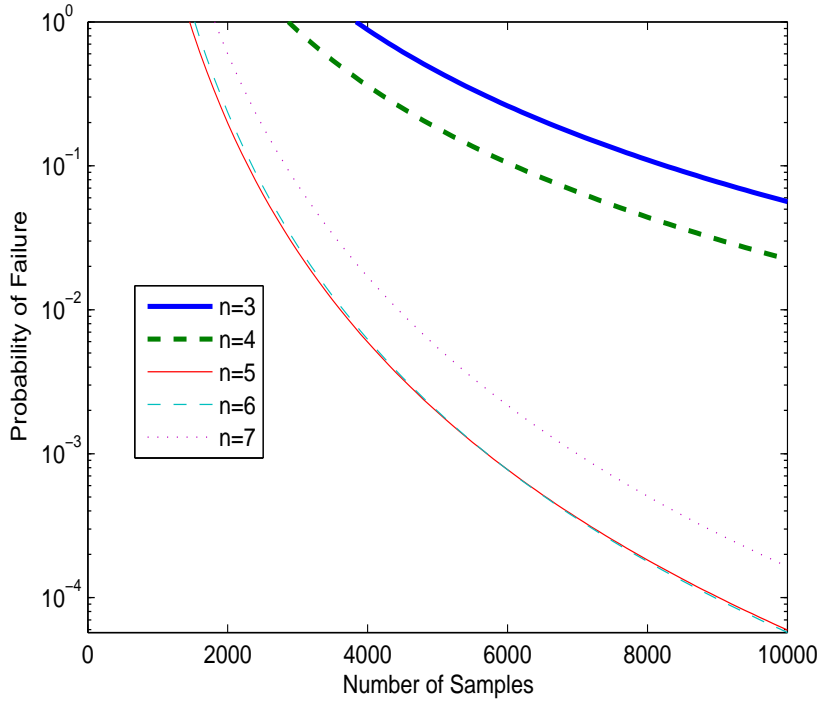}}
\caption{
Bounds for the probability of failure of exact reconstruction due to
Theorem \ref{thm_main0} with $M=10$, $D=10000$.}
\end{center}
\end{figure}

In figure \ref{fig_prob_thm0} we have chosen 
$M=10$, $D=10000$ and $n=3,4,5,6,7$ to show a logarithmic plot of the probability 
bound (\ref{Pbound0}) of Theorem \ref{thm_main0} 
versus the number of samples. The parameter $\beta$ was chosen always
near to $1/2$ and then $\kappa$ was determined such that there is equality
in (\ref{cond_kappa_0}). One can see clearly, that here 
$n=5$ or $n=6$ is the optimal choice depending on the precise value of the
number of samples $N$. Unfortunately, it seems that these bounds are
quite pessimistic when compared to the 
numerical experiments (see next section). In the given example 
one needs at least about $N=2000$ samples 
(corresponding to a ``non-linear oversampling factor'' of $200$)
in order that the bound becomes non-trivial.    

\begin{figure}\label{fig_prob_estim}
\parbox[t]{16cm}{
\includegraphics[width=8cm,height=6cm]{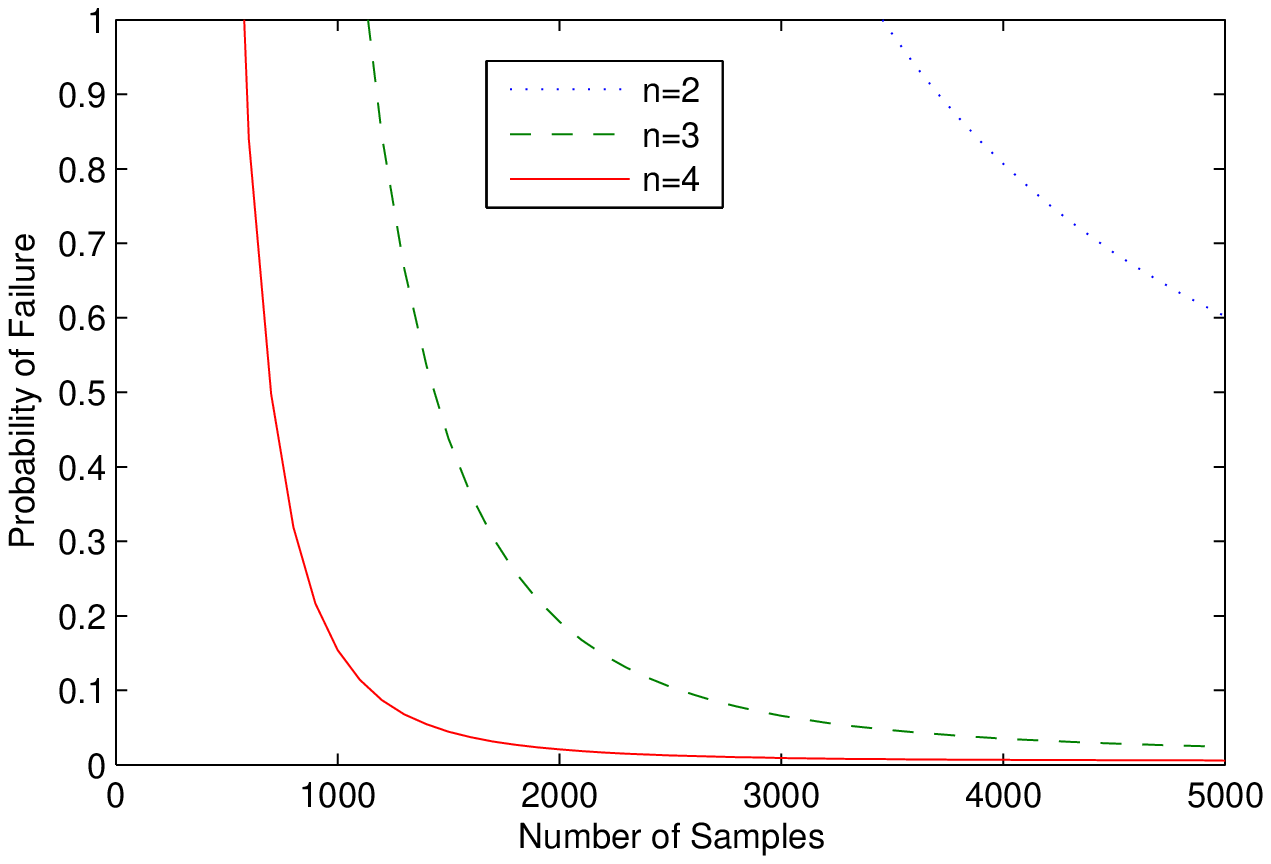}
\includegraphics[width=8cm,height=6cm]{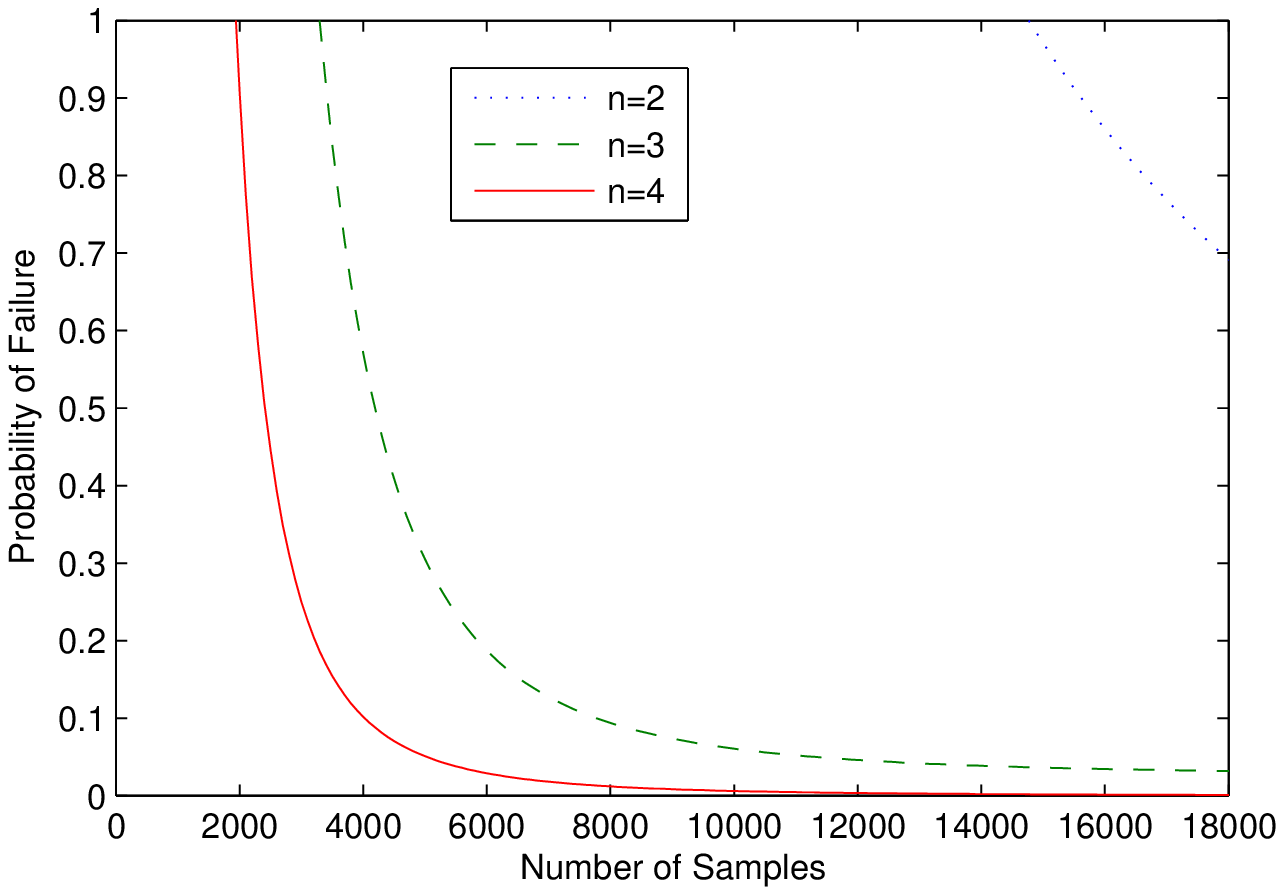}
}
\caption{Probability of failure of exact reconstruction
for $\E|T| = 4$, $D=5000$ (left) and $\E|T| = 8, D=20000$ (right)
due to Theorem \ref{thm_main}}
\end{figure}

Based on the computation of the explicit values of the numbers $Q$ and
$Q^*$ we can also illustrate the probability bound (\ref{Pbound})
in Theorem \ref{thm_main}. Unfortunately, we may
only take $n\leq 4$ since for higher values of $n$ the corresponding
numbers $Q$ and $Q^*$ are not at our disposal.
Figure (\ref{fig_prob_estim}) shows a plot of the bound
(\ref{Pbound}). We have chosen $n=2,3,4$ and 
$(\E|T|,D) = (4,5000),\, (8,20000)$
and varied the number $N$ of sampling points. For $n=2$ we have chosen
$K_1=2, K_2=1$, for $n=3$: $K_1=3,K_2=2,K_3=1$ and for $n=4$ we took
$K_1=4,K_2=2, K_3=1,K_4=1$ as suggested in Remark \ref{rem_thm}(a).
It turned out that good choices for $\beta$ are around $1/2$ and for
$\kappa\approx 10^{-3}$ (with slight variations for the different
choices of the other parameters). The remaining parameter $\alpha$ was chosen
such that there is equality in (\ref{cond_kappa_thm}). 

\begin{figure}\label{fig_compare}
\begin{center}
\parbox[t]{12cm}{
\includegraphics[width=12cm,height=8cm]{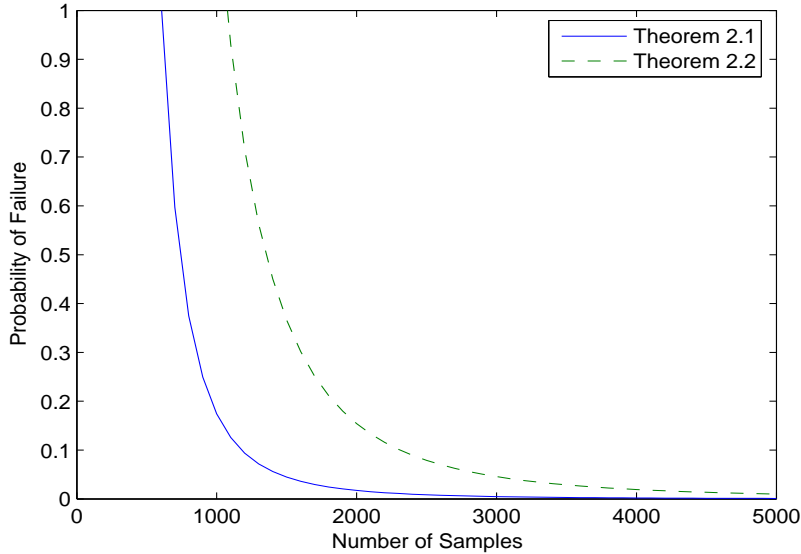}
}
\caption{Comparison of the bounds of Theorem \ref{thm_main0} and Theorem 
\ref{thm_main} for $M=\E|T|=4$, $D=5000$ and $n=4$.}
\end{center}
\end{figure}

Looking at the plot one 
realizes clearly that the bound becomes better for larger $n$. 
However, as above 
the bounds are still quite pessimistic. Nevertheless, as already remarked 
one expects them to be at least better than the ones of 
Theorem \ref{thm_main0}.
Figure \ref{fig_compare} supports this intuition. Indeed, we plotted
the different bounds for $M=\E|T|=4$, 
$D=5000$, $n=4$ and $K_1=4,K_2=2,K_3=1$ and
$K_4=1$. Apparently the curve for the bound of Theorem \ref{thm_main}
is far below the one of Theorem \ref{thm_main0}. 
Unfortunately, we cannot yet use the full strength of Theorem \ref{thm_main}
as we are still lacking an efficient way to actually compute
the bound explicitly for higher values of $n$. Actually up to now Theorem
\ref{thm_main0} still gives the better bound in most situations because
we are able to evaluate (\ref{Pbound0}) for arbitrary $n$.

Let us finally discuss possible reasons why 
the theoretical bounds are quite pessimistic.
Both theorems 
give bounds for the probability that exact reconstruction
holds for {\it all} choices of the coefficients $f$ on $T$, 
while the numerical experiments
in the next section choose also the coefficients on $T$ at random. 
(Of course, it is impossible to check all possible coefficients 
by some algorithm.) 
Intuitively, it is very plausible that in such an experiment 
the probability of failure of exact
reconstruction is much lower than for the situation
in our main Theorem \ref{thm_main}. We remark that it 
seems to be an interesting project to investigate theoretically 
also the case that the coefficients of $f$ on $T$ are chosen at random, 
see also Section 5 in \cite{CT}. We plan to pursue this issue
in a follow-up paper.

Of course, the theoretical bounds may also be pessimistic compared
to reality since some of the estimates in the proof are perhaps not sharp.
However, it seems to be hard to improve on the method of our proof.

\section{Numerical experiments}
\label{sec_numerical}

Let us describe some numerical tests of the proposed 
sampling resp. reconstruction method. In order to use convex 
optimization techniques
we reformulate the optimization problem \ref{min_problem} as
the following equivalent problem,
\begin{align}
\label{mod_problem}
\min \sum_k u_k \quad \mbox{ subject to } & \sqrt{(c^{(1)}_k)^2 + (c^{(2)}_k)^2} \leq u_k,\\
& \sum_{k} (c_k^{(1)} + i c_k^{(2)})e^{ik\cdot x_j} = f(x_j)\notag
\end{align}
with $u_k$ and $c^{(1)}_k$ and $c^{(2)}_k$, $k \in [-q,q]^d$, as
real optimization variables. The solution to the original problem
\ref{min_problem} is then given as $c_k = c^{(1)}_k + i c_k^{(2)}$.  

A problem of the above type (\ref{mod_problem}) is known
as second order cone program \cite{BV}. Efficient algorithms to
solve such problems exist. We have 
used the toolbox MOSEK (in connection with MATLAB), 
which provides an interior point solver
for cone problems.
We remark that if the coefficients $c_k$ 
are real-valued 
then the minimization problem (\ref{min_problem}) can be recast as a linear
program. 

Our numerical experiment has the following form. We first choose
the sparsity $M$, the maximal degree $q$ (we only tested for $d=1$)
and the number of samples $N$. Then the following steps are done:
\begin{enumerate}\itemsep-0.5pt
\item Choose a random subset $T \subset [-q,q]$ of size $M$ from the 
uniform distribution. (Generate a random permutation of $[-q,q]$ and take
the first $M$ elements.)
\item Randomly generate the coefficients $c_k$ for $k \in T$ by choosing
their real part and imaginary part from
a standard normal distribution.
\item Randomly select $x_1,\hdots,x_N$ independently 
from the uniform distribution on
$[0,2\pi]$.
\item Generate $f(x_j) = \sum_{k \in T} c_k e^{ikx_j}$, 
$j=1,\hdots,N$.
\item Solve the minimization problem (\ref{mod_problem}).
\item Compare the result to the original vector of coefficients.
\end{enumerate}

\begin{figure}\label{fig_numerical}
\center{
\parbox[t]{12cm}{
\includegraphics[width=12cm]{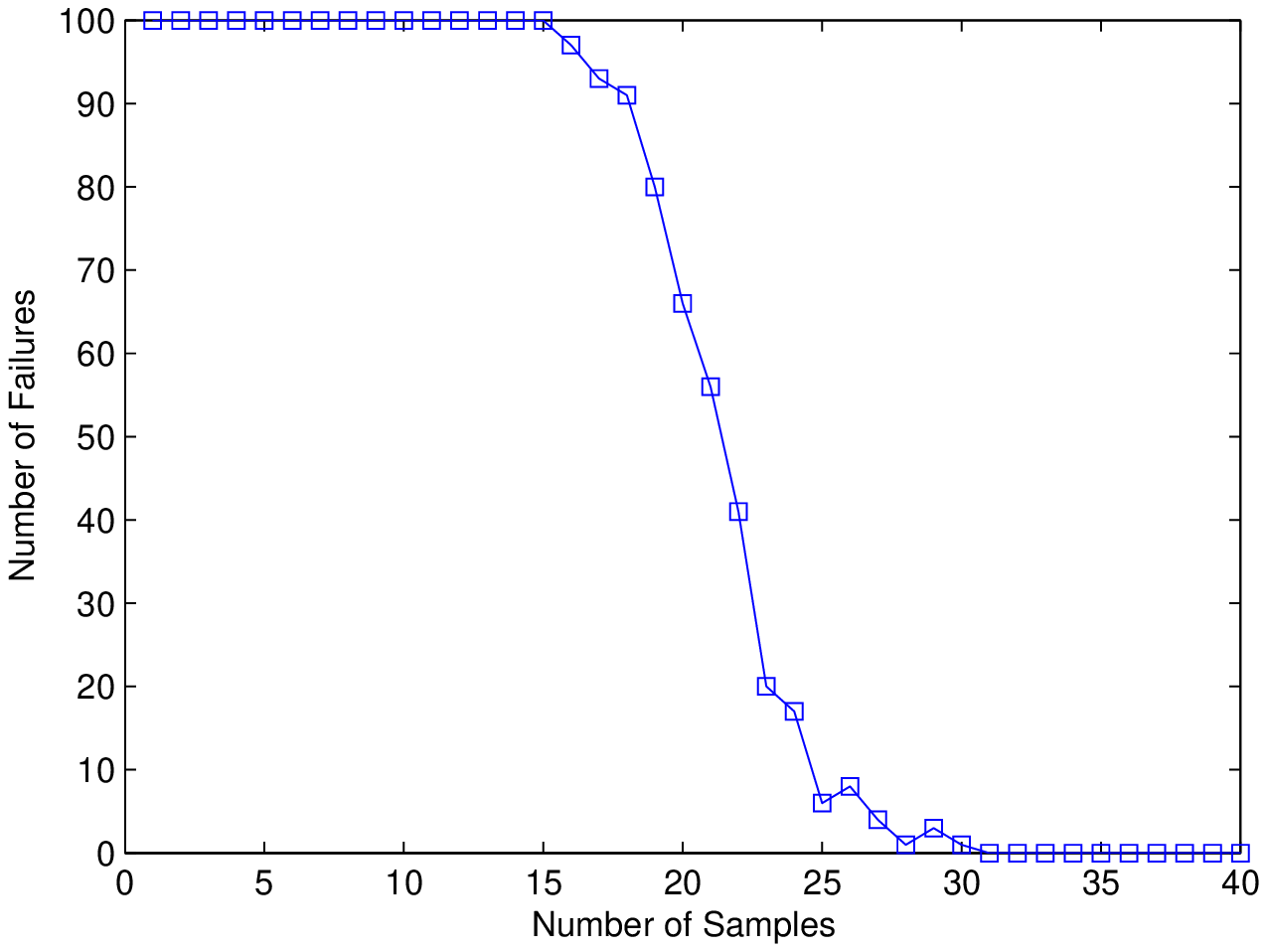}
}}
\center{Numerical results: number of failures out of $100$ 
trials for $M=|T|=8$ and $q=40$ versus number of samples}
\end{figure}

For figure \ref{fig_numerical} we have 
chosen $q=40$, i.e., $D= (2q+1) = 81$ and $M=|T|=8$. 
Then for each $N$ between $1$ and $40$ 
we ran the above procedure
$100$ times and counted how often exact reconstruction failed. The result
is illustrated in the plot. As one can see for $N$ larger than $30$ 
(corresponding to a non-linear oversampling factor of about $4$)
our reconstruction method always succeeded in giving back
the original function exactly!

Comparing these results with the bounds of Theorem \ref{thm_main}
as illustrated in the previous section one realizes that in practice the 
method works even much better than we are able to predict theoretically.
So this method seems to have quite a lot of potential for practical
applications of signal reconstruction. 

\begin{appendix}
\section{Appendix}
\subsection{Tables for $Q(n,t,s,R)$ and $Q^*(K,m,t,s,R)$}

In the following tables we list some values for the numbers 
$Q(n,t,s,R)$ and $Q^*(K,m,t,s,R)$ that were computed by the procedures
described in Section \ref{sec_Qn}. For the probability estimation 
(\ref{Pbound}) the numbers $Q^*(2K_m,m,t,s,R)$, $m=1,\hdots,n$ 
are needed and we have chosen $n=4$ and $K_1=4, K_2 = 2, K_3=1, K_4=1$ 
(since for higher numbers of $n$ and $K_m$ computing 
times are absurdly long).
So the numbers $Q^*(K,m,t,s,R)$ have to be computed 
for $(K,m) = (8,1), (4,2), (2,3), (2,4)$.

Note that $Q(n,t,s,R) = 0$ if $R \geq \min\{t,s\}$ or $s = 1$ or 
$t > n/2$ and $Q^*(n,t,s,R) = 0$ if $R > \min \{t,s\}$, 
which is the reason why we do not reproduce these cases. Also 
recall that $U(n,1) = \emptyset$ and, unless $m=1$, $U(K,m,1) = \emptyset$,
hence, $Q(n,t,1,R) = 0$ and $Q^*(K,m,t,1,R) = 0$.

\bigskip

\noindent
\begin{tabular}{|l||c|} \hline
$Q(4,1,s,R)$ & $R = 0$ \\ \hline\hline
$s = 2$      &   1  \\ \hline  
$s = 3$      &   2  \\ \hline
$s = 4$      &   1  \\ \hline
\end{tabular}
\begin{tabular}{|l||c|c|} \hline
$Q(4,2,s,R)$ & $R = 0$ & $R=1$ \\ \hline\hline
$s = 2$      &   2     &   1   \\ \hline  
$s = 3$      &   2     &   4   \\ \hline
$s = 4$      &   0     &   3   \\ \hline
\end{tabular}

\bigskip

\noindent
\begin{tabular}{|l||c|} \hline
$Q(6,1,s,R)$ & $R = 0$  \\ \hline\hline
$s = 2$      &    1     \\ \hline  
$s = 3$      &    10     \\ \hline
$s = 4$      &    20     \\ \hline
$s = 5$      &    9      \\ \hline 
$s = 6$      &    1      \\ \hline
\end{tabular}
\begin{tabular}{|l||c|c|} \hline
$Q(6,2,s,R)$ & $R = 0$ & $R=1$ \\ \hline\hline
$s = 2$      &    9    &   16     \\ \hline  
$s = 3$      &    46    &  204    \\ \hline
$s = 4$      &    45    &  455    \\ \hline
$s = 5$      &    9     &  216     \\ \hline 
$s = 6$      &    0     &  25      \\ \hline
\end{tabular}
\begin{tabular}{|l||c|c|c|} \hline
$Q(6,3,s,R)$ & $R = 0$ & $R=1$ & $R=2$ \\ \hline\hline
$s = 2$      &    6    &   9   &   0    \\ \hline  
$s = 3$      &    15    &  78  &   57    \\ \hline
$s = 4$      &    5    &   87  &   208    \\ \hline
$s = 5$      &    0     &  18  &   117    \\ \hline 
$s = 6$      &    0     &  0   &   15    \\ \hline
\end{tabular}

\bigskip

\noindent
\begin{tabular}{|l||c|} \hline
$Q(8,1,s,R)$ & $R = 0$  \\ \hline\hline
$s = 2$      &    1    \\ \hline  
$s = 3$      &    42     \\ \hline
$s = 4$      &    231     \\ \hline
$s = 5$      &    294      \\ \hline 
$s = 6$      &    126      \\ \hline
$s = 7$      &    20    \\\hline
$s = 8$      &    1     \\\hline
\end{tabular}
\begin{tabular}{|l||c|c|} \hline
$Q(8,2,s,R)$ & $R = 0$ & $R=1$ \\ \hline\hline
$s = 2$      &    34    &  85     \\ \hline  
$s = 3$      &    674   &  4324   \\ \hline
$s = 4$      &    1970  &  25519  \\ \hline
$s = 5$      &    1386  &  33600  \\ \hline 
$s = 6$      &    308   &  14686  \\ \hline
$s = 7$      &    20    &  2360   \\ \hline
$s = 8$      &    0     &  119    \\ \hline 
\end{tabular}
\begin{tabular}{|l||c|c|c|} \hline
$Q(8,3,s,R)$ & $R = 0$ & $R=1$ & $R=2$ \\ \hline\hline
$s = 2$      &    72    &  418    &   0    \\ \hline  
$s = 3$      &    732   &  892    &   10896   \\ \hline
$s = 4$      &    1218  &  27446  &   84526  \\ \hline
$s = 5$      &    504   &  19944  &   123612  \\ \hline 
$s = 6$      &    56    &  4556   &   57128   \\ \hline
$s = 7$      &    0     &   304   &   9496    \\ \hline
$s = 8$      &    0     &    0    &   490     \\ \hline
\end{tabular}
\begin{tabular}{|l||c|c|c|c|} \hline
$Q(8,4,s,R)$ & $R = 0$ & $R=1$   & $R=2$   & $R=3$ \\ \hline\hline
$s = 2$      &    24   &  81     &   0     &   0   \\ \hline  
$s = 3$      &    112  &  1208   &   3090  &   0   \\ \hline
$s = 4$      &    84   &  2018   &   11944 &   10209    \\ \hline
$s = 5$      &    14   &  800    &   9368  &   20688\\ \hline 
$s = 6$      &    0    &  86     &   2236  &   10908 \\ \hline
$s = 7$      &    0    &   0     &   156   &   1944  \\ \hline
$s = 8$      &    0    &   0     &    0    &   105   \\ \hline
\end{tabular}

\bigskip

\newpage

\noindent
\underline{$Q^*(8,1,t,s,R)$:}\\
\begin{tabular}{|l||c|c|} \hline
$t=1$            & $R=0$ & $R=1$ \\ \hline\hline
$s=1$            &  1    & 0     \\ \hline
$s=2$            &  34   & 93    \\ \hline
$s=3$            &  72   & 894   \\ \hline
$s=4$            &  24   & 1677  \\ \hline
$s=5$            &  0    & 1050  \\ \hline
$s=6$            &  0    & 266   \\ \hline
$s=7$            &  0    & 28    \\ \hline
$s=8$            &  0    & 1     \\ \hline
\end{tabular}
\begin{tabular}{|l||c|c|c|} \hline
$t=2$ & $R=0$ & $R=1$ & $R=2$ \\ \hline\hline
$s=1$ &   34  &   85  &  0    \\ \hline
$s=2$ &  610  & 6598  & 7905  \\ \hline
$s=3$ &  792  & 20420 & 93742 \\ \hline
$s=4$ &  168  & 13736 & 188515 \\ \hline
$s=5$ &  0    & 3380  & 121570 \\ \hline
$s=6$ &  0    & 408   & 31246  \\ \hline
$s=7$ &  0    & 16    & 3316   \\ \hline
$s=8$ &  0    & 0     & 119    \\ \hline
\end{tabular}\\
\begin{tabular}{|l||c|c|c|c|} \hline
$t=3$ & $R=0$ & $R=1$ & $R=2$ & $R=3$ \\ \hline\hline
$s=1$ &  72   &  418  &   0   &  0    \\ \hline
$s=2$ & 792   & 14572 & 46866 &  0    \\ \hline
$s=3$ & 792   & 25704 & 210638& 236206 \\ \hline 
$s=4$ & 144   & 10104 & 192512& 630730 \\ \hline
$s=5$ &  0    & 1368  & 63134 & 449998 \\ \hline
$s=6$ &  0    &  72   & 8700  & 121568 \\ \hline
$s=7$ &  0    &  0    & 400   & 13320  \\ \hline
$s=8$ &  0    &  0    &  0    & 490    \\ \hline
\end{tabular}
\begin{tabular}{|l||c|c|c|c|c|} \hline
$t=4$ & $R=0$ & $R=1$ & $R=2$ & $R=3$ & $R=4$ \\ \hline\hline
$s=1$ &   24  &  81   &  0    &   0   &   0   \\ \hline
$s=2$ &  168  & 3744  &  9423 &  0    &   0   \\ \hline
$s=3$ &  144  & 4440  &  38472&  58374&   0   \\ \hline
$s=4$ &   24  & 1296  &  21060&  96384&   59841 \\ \hline
$s=5$ &   0   & 96    &  3816 &  37302&  69036  \\ \hline
$s=6$ &   0   &  0    &  216  & 5292  &  22422  \\ \hline
$s=7$ &   0   &  0    &  0    & 240   &  2700   \\ \hline
$s=8$ &   0   &  0    &   0   &  0    &  105    \\ \hline
\end{tabular}

\bigskip 

\noindent
\underline{$Q^*(4,2,t,s,R)$:}\\
\begin{tabular}{|l||c|c|} \hline
$t=1$ & $R=0$ & $R=1$ \\ \hline\hline
$s=2$            &   3   &  5    \\ \hline
$s=3$            &   37  &  171  \\ \hline
$s=4$            &   56  &  596  \\ \hline 
$s=5$            &   21  &  555  \\ \hline
$s=6$            &   2   &  186  \\ \hline
$s=7$            &   0   &   24  \\ \hline
$s=8$            &   0   &   1   \\ \hline
\end{tabular}
\begin{tabular}{|l||c|c|c|} \hline
$t=2$ & $R=0$ & $R=1$ & $R=2$ \\ \hline\hline
$s=2$            &   47  & 415   & 490   \\ \hline
$s=3$            &   274 & 5866  & 18612 \\ \hline 
$s=4$            &   226 & 9537  & 67825 \\ \hline
$s=5$            &   46  & 3946  & 64552 \\ \hline
$s=6$            &   2   & 480   & 21890 \\ \hline
$s=7$            &   0   & 12    & 2844  \\ \hline
$s=8$            &   0   & 0     & 119   \\ \hline
\end{tabular}\\
\begin{tabular}{|l||c|c|c|c|} \hline
$t=3$ & $R=0$ & $R=1$ & $R=2$ & $R=3$ \\ \hline\hline
$s=2$            &  50   &  744  &  3126  &  0   \\ \hline
$s=3$            &  134  & 4930  &  42410 &  54446  \\ \hline
$s=4$            &  54   & 4070  &  70998 &  244358 \\ \hline
$s=5$            &   4   & 776   &  31452 &  250008 \\ \hline
$s=6$            &  0    & 26    &  4422  &  87672  \\ \hline
$s=7$            &  0    & 0     &  166   &  11594  \\ \hline
$s=8$            &  0    & 0     &  0     &  490    \\ \hline
\end{tabular}
\begin{tabular}{|l||c|c|c|c|c|} \hline
$t=4$ & $R=0$ & $R=1$ & $R=2$ & $R=3$ & $R=4$ \\ \hline\hline
$s=2$            &  8    &  108  &   724 &  0    &  0    \\ \hline
$s=3$            &  10   &  412  &  4910 & 16508 &  0    \\ \hline
$s=4$            &  2    &  186  &  4377 & 30778 &  33117\\ \hline
$s=5$            &  0    &  17   &  962  & 14607 & 44894 \\ \hline
$s=6$            &  0    &  0    &  53   & 2266  & 17421 \\ \hline
$s=7$            &  0    &  0    &  0    & 102   & 2418  \\ \hline
$s=8$            &  0    &  0    &  0    &  0    & 105   \\ \hline
\end{tabular}
\bigskip

\newpage

\noindent
\underline{$Q^*(2,3,t,s,R)$:}\\
\begin{tabular}{|l||c|c|} \hline
$t=1$ & $R=0$ & $R=1$ \\ \hline\hline
$s=2$ &   1   &  1   \\ \hline        
$s=3$ &   7   &  15  \\ \hline
$s=4$ &   6   &  25  \\ \hline
$s=5$ &   1   &  10  \\ \hline
$s=6$ &   0   &  1   \\ \hline
\end{tabular}
\begin{tabular}{|l||c|c|c|} \hline
$t=2$ & $R=0$ & $R=1$ & $R=2$ \\ \hline\hline
$s=2$ &  3    &  27   &  20   \\ \hline
$s=3$ & 10    & 200   & 340   \\ \hline
$s=4$ & 4     & 172   & 599   \\ \hline
$s=5$ & 0     & 29    & 246   \\ \hline
$s=6$ & 0     & 0     & 25    \\ \hline
\end{tabular}
\begin{tabular}{|l||c|c|c|c|} \hline
$t=3$ & $R=0$ & $R=1$ & $R=2$ & $R=3$ \\ \hline\hline
$s=2$ &  1    &   9   &  20   &  0    \\ \hline
$s=3$ &  1    &   23  & 166   &  140  \\ \hline
$s=4$ &  0    &   7   & 132   &  326  \\ \hline
$s=5$ &  0    &   0   & 21    &  144  \\ \hline
$s=6$ &  0    &   0   &  0    &  15   \\ \hline
\end{tabular}
\bigskip

\noindent
\underline{$Q^*(2,4,t,s,R)$:}\\
\begin{tabular}{|l||c|c|} \hline
$t=1$ & $R=0$ & $R=1$ \\ \hline\hline
$s=2$ &  1    &  1  \\ \hline
$s=3$ &  31   &  63 \\ \hline
$s=4$ & 90    &  301 \\ \hline
$s=5$ & 65    &  350  \\ \hline
$s=6$ & 15    &  140  \\ \hline
$s=7$ & 1     &  21   \\ \hline
$s=8$ & 0     &  1    \\ \hline
\end{tabular}
\begin{tabular}{|l||c|c|c|} \hline
$t=2$ & $R=0$ & $R=1$ & $R=2$ \\ \hline\hline
$s=2$ &  11   &  135  &  92   \\ \hline
$s=3$ &  157  & 4225  &  6804 \\ \hline
$s=4$ &  222  & 11981 &  34326 \\ \hline
$s=5$ &  69   & 8438  &  40878 \\ \hline
$s=6$ &  5    & 1902  &  16538 \\ \hline
$s=7$ &  0    & 124   &  2494  \\ \hline
$s=8$ &  0    &  0    &  119   \\ \hline
\end{tabular}\\
\begin{tabular}{|l||c|c|c|c|} \hline
$t=3$ & $R=0$ & $R=1$ & $R=2$ & $R=3$ \\ \hline\hline
$s=2$ &  11   &  223  &  746  & 0     \\ \hline
$s=3$ &  71   &  2974 &  23519& 19496 \\ \hline
$s=4$ &  48   &  3907 & 63319 & 124316 \\ \hline
$s=5$ &  5    &  1171 & 42404 & 159770 \\ \hline
$s=6$ &  0    &  85   & 9135  & 66730  \\ \hline
$s=7$ &  0    &  0    & 572   & 10208  \\ \hline
$s=8$ &  0    &  0    &  0    & 490    \\ \hline
\end{tabular}
\begin{tabular}{|l||c|c|c|c|c|} \hline
$t=4$ & $R=0$ & $R=1$ & $R=2$ & $R=3$ & $R=4$ \\ \hline\hline
$s=2$ &  2    &  36   &  172  &  0    &   0   \\ \hline
$s=3$ &  4    &  164  &  2393 &  7309 &   0   \\ \hline
$s=4$ &  1    &  101  &  2865 &  20438& 17650 \\ \hline 
$s=5$ &  0    &  11   &  820  &  12988& 29756 \\ \hline
$s=6$ &  0    &  0    &  58   &  2643 & 13574 \\ \hline
$s=7$ &  0    &  0    &  0    &  156  & 2154  \\ \hline
$s=8$ &  0    &  0    &  0    &  0    &  105  \\ \hline
\end{tabular}

\end{appendix}

\end{document}